\documentclass[12pt]{article}
\usepackage{amsfonts}

\newcommand{\calC}{\mathcal{C}}
\newcommand{\calV}{\mathcal{V}}
\newcommand{\calD}{\mathcal{D}}
\newcommand{\calM}{\mathcal{M}}
\newcommand{\calE}{\mathcal{E}}
\newcommand{\calU}{\mathcal{U}}
\newcommand{\calK}{\mathcal{K}}
\newcommand{\calS}{\mathcal{S}}
\newcommand{\calB}{\mathcal{B}}
\newcommand{\calL}{\mathcal{L}}
\newcommand{\boldd}{{\rm\bf d}}
\newcommand{\boldx}{{\rm\bf x}}
\newcommand{\boldz}{{\rm\bf z}}
\newcommand{\bolde}{{\rm\bf e}}
\newcommand{\dis}{\displaystyle}
\newcommand{\rarrow}{\rightarrow}
\renewcommand{\theequation}{\thesection.\arabic{equation}}

\newtheorem{theorem}{Theorem}[section]
\newtheorem{lemma}[theorem]{Lemma}

\newtheorem{corollary}[theorem]{Corollary}
\newtheorem{proposition}[theorem]{Proposition}
\newtheorem{remark}[theorem]{Remark}

\begin{document}

\title{Linear Relations in the Calkin Algebra \\ for Composition Operators}

\author{Thomas Kriete\thanks{Work of the first author was supported in part by a Sesquicentennial Associateship at the University of Virginia.}
 \\[2mm]
{\small Department of Mathematics} \\[-1mm] {\small University of Virginia} \\[-1mm]
{\small Charlottesville, VA  22904}
\and
Jennifer Moorhouse \\[2mm] {\small Department of Mathematics} \\[-1mm] {\small Colgate University} \\[-1mm] {\small Hamilton, NY 13356}}

\date{}
\maketitle

\begin{abstract}
We consider this and related questions: When is a finite linear combination of composition operators, acting on 
the Hardy space or the standard weighted Bergman spaces on the unit disk, a compact operator?
\end{abstract}

\maketitle

\section{Introduction}
\renewcommand{\theequation}{\thesection.\arabic{equation}}

For $\beta \geq 1$,
let $\calD_\beta$ denote the reproducing kernel Hilbert space of functions analytic in the unit disk
$D = \{z: |z| < 1\}$ and having the kernel functions $k_w(z) = (1-\overline{w}z )^{-\beta}$. Thus, $f(w) = \langle f,k_w\rangle$ for $w$ in $D$ and $f$ in $\calD_\beta$.  The Hardy space $H^2$ is exactly $\calD_1$ and when $\beta > 1$, $\calD_\beta$ is the standard weighted Bergman space $A^2_\alpha$ with $\alpha+2 = \beta$, see Section 2.1. We consider composition operators $C_\varphi: f \rarrow f \circ \varphi$ acting on $\calD_\beta$, where $\varphi$ is an analytic self-map of $D$. When $\beta \geq 1$, every $C_\varphi$ lies in $\calB(\calD_\beta)$, the algebra of bounded operators on $\calD_\beta$. Unlike the classes of Toeplitz and Hankel operators which act on $\calD_\beta$, the set of composition operators in $\calB(\calD_\beta)$ has no obvious additive or linear structure. However, in the Bergman space case $\beta > 1$, the second author observed additive structure modulo the ideal $\calK$ of compact operators and characterized those pairs $\varphi$ and $\psi$ for which $C_\varphi - C_\psi$ is compact \cite{Mo}. Our purpose here is twofold: to present some analogous results for the $H^2$ case $\beta =1$, and to pass from additive to linear structure in the Calkin algebra $\calB(\calD_\beta)/\calK$.

Composition operators which are themselves compact were characterized  in the $A^2_\alpha$ case by MacCluer and Shapiro \cite{MS} and on $H^2$ by Shapiro
\cite{S}; in \cite{Sa1} Sarason found a different 
condition, sufficient for $H^2$ and necessary and sufficient for $L^1$,  later shown by Shapiro and Sundberg \cite{SSu2} to be necessary in the $H^2$ case as well.
The problem of compact difference was raised in explicit form by Shapiro and Sundberg \cite{SSu1} and MacCluer \cite{M}; these authors found several criteria, some necessary and some sufficient. More recently, MacCluer, Ohno, and Zhao \cite{MOZ} have shown that compactness of $C_\varphi-C_\psi$ acting on $H^\infty$, the space of bounded analytic functions on $D$, is characterized in terms of the quantity $\rho(z) = \left|\frac{\varphi(z)-\psi(z)}{1-\overline{\varphi(z)} \psi(z)}\right|$, the pseudo-hyperbolic distance between values of $\varphi$ and $\psi$. The second author of the present article subsequently found that $\rho$ plays a related role on the Bergman spaces, as follows:

\renewcommand{\thetheorem}{A}
\begin{theorem} {\rm \cite{Mo}} Suppose that $\varphi$ and $\psi$ are analytic self-maps of $D$ and $\alpha > -1$. Then $C_\varphi-C_\psi$ is a compact operator on $A^2_\alpha$ if and only if
$$
\lim_{|z| \rarrow 1} \rho(z) \left\{\frac{1-|z|^2}{1-|\varphi(z)|^2} + 
\frac{1-|z|^2}{1-|\psi(z)|^2}\right\} = 0.
$$
\end{theorem}

\noindent
Let $F(\varphi)$ denote the set of points in the unit circle $\partial D$ at which $\varphi$ has a finite angular derivative in the sense of Caratheodory, see Section 2.2. Also, we will use the notation $A \equiv B$ (mod $\calK$) to indicate that two bounded operators $A$ and $B$ have compact difference. There is a ``sum theorem,'' as follows:

\renewcommand{\thetheorem}{B}
\begin{theorem} {\rm \cite{Mo}}
Let $\varphi, \varphi_1,\ldots,\varphi_n$ be analytic self-maps of $D$ for which the sets $F(\varphi_i)$, $i = 1,\ldots,n$, are pairwise disjoint and with $F(\varphi) = F(\varphi_1) \cup \cdots \cup F(\varphi_n)$. Consider $C_\varphi,C_{\varphi_{1}},\ldots,C_{\varphi_{n}}$ as acting on $A^2_\alpha$ where $\alpha > -1$, and let
$$
\rho_i = \left|\frac{\varphi-\varphi_i}{1-\overline{\varphi} \varphi_i}\right|, \quad i = 1,\ldots,n.
$$
Then the following are equivalent:

(i) For each $i = 1,\ldots,n$ and each $\zeta$ in $F(\varphi_i)$,
$$
\lim_{z \rarrow \zeta} \rho_i(z) \left\{\frac{1-|z|^2}{1-|\varphi(z)|^2} + 
\frac{1-|z|^2}{1-|\varphi_i(z)|^2}\right\} = 0.
$$

(ii) $C_\varphi \equiv C_{\varphi_{1}} + \cdots + C_{\varphi_{n}}
\ (\mbox{\rm mod }\calK).$
\end{theorem}

\noindent
The following useful modification of Theorem A, which localizes the notion of compact difference, is implicit in \cite{Mo}.

\renewcommand{\thetheorem}{C}
\begin{theorem} 
Suppose that $\varphi$ and $\psi$ are analytic self-maps of $D$, $\alpha > -1$ and $G$ is a measurable subset of $D$. If 
$$
\lim_{|z| \rarrow 1} \chi_{_{G}} (z) \rho(z) \left\{\frac{1-|z|^2}{1-|\varphi(z)|^2} + \frac{1-|z|^2}{1-|\psi(z)|^2}\right\} = 0,
$$
then $M_{\chi_{_{G}}}(C_\varphi-C_\psi)$ is a compact operator from $A^2_\alpha$ into its containing $L^2$ space. Here  $M_{\chi_{_{G}}}$ denotes the operator of multiplication by the characteristic function $\chi_{_{G}}$.
\end{theorem}

Other recent related work includes the following:  Bourdon, Levi, Narayan, and J.~H. Shapiro \cite{BLNS} show that a composition operator associated with an ``almost linear fractional'' map is, in fact, a compact perturbation of a linear fractional composition operator; Bourdon \cite{B} treats the question of compact difference vs.~topological connectedness  for linear fractional maps (see Section 6 below); J.~E.~Shapiro
\cite{Sh} shows that if $C_\varphi-C_\psi$ is compact on $H^2$, the singular parts of the Clark measures of $\varphi$ and $\psi$ coincide
(see Sect.~2.3).  Most recently, Nieminen and Saksman \cite{NiSak} have shown that the just-quoted condition of J.~E.~Shapiro on singular parts, plus uniform integrability of the differences of the absolutely continuous densities of the respective Clark measures of $\varphi$ and $\psi$, are together equivalent to compactness of $C_{\varphi}-C_\psi$ on $L^1$ (or, on the space of complex Borel measures on $\partial\mathbb{D}$) and are thus sufficient on $H^2$. In a different direction, Gorkin and Mortini \cite{GMor} have characterized compactness for finite linear combinations of composition operators acting on  uniform algebras.

Our analogues of Theorems A, B, and C for $H^2$ require some different methods.
For us the key is an application of Clark measures, following ideas of Sarason \cite{Sa1}, Cima and Matheson \cite{CiMa} and J.~E. Shapiro \cite{Sh}. In Section 3 we obtain essential norm estimates for weighted composition operators on $H^2$ analogous to the Cima-Matheson essential norm formula for (unweighted) composition operators \cite{CiMa}. We combine these results in Section 4 with a general principle of Moorhouse and Toews \cite{MoT} and Carleson measure estimates as in \cite{Mo} to obtain our $H^2$ results.

Section 5 is devoted to the question of when a given finite linear combination of composition operators is compact. We look at lower bounds, given in terms of first- and higher-order boundary data, for the essential norm of a linear combination; these results further develop ideas of MacCluer \cite{M}. 
We introduce the class $\calS$ of analytic self-maps $\varphi$ of $D$ having ``sufficient data'' at every point in $\partial D$ where $\varphi$ makes significant contact with the boundary. For $\varphi$ and $\psi$ in $\calS$ we make clear the obstructions to the essential norm $\|C_\varphi - C_\psi\|_e$ being small, the conditions under which it must be small, and when it is zero, that is, when $C_\varphi-C_\psi$ is compact.
For $\varphi_1,\ldots,\varphi_n$ in $\calS$, we characterize, via a finite system of linear equations involving boundary data of these maps, those coefficients for which $c_1 C_{\varphi_{1}} + \cdots + c_n C_{\varphi_{n}}$ is compact. An application is a simple algorithm for determining the dimension of the vector space in $\calB(\calD_\beta)/\calK$ spanned by the cosets $[C_{\varphi_{1}}],\ldots,[C_{\varphi_{n}}]$.

The final Section 6 concerns a problem first studied by Berkson \cite{Be}, subsequently considered by MacCluer \cite{M}, Shapiro and Sundberg \cite{SSu1}, and most recently by Toews and the second author \cite{MoT} and Bourdon \cite{B}: to characterize those pairs $\varphi$ and $\psi$ for which $C_\varphi$ and $C_\psi$ lie in the same connected component of the topological space of composition operators, equipped with the norm topology on $\calB(\calD_\beta)$. We observe that a general sufficient condition of the second author for the Bergman space case \cite{Mo} extends to $H^2$ and apply this result to those $C_\varphi$ with $\varphi$ lying in a certain subclass $\calS_0$ of $\calS$ to determine when $C_\varphi$ and $C_\psi$ lie in the same component of $\{C_\varphi: \varphi \in \calS_0\}$.

A variation on ideas of Berkson \cite{Be}, Shapiro and Sundberg \cite{SSu1}, and MacCluer \cite{M} (see Exercise 9.3.2 in \cite{CM})
states that if $\varphi_1,\ldots,\varphi_n$ are analytic self-maps of $D$, and if $J(\varphi)$ denotes the set of points $e^{i\theta}$ in $\partial D$ with $|\varphi(e^{i\theta})| =1$, then
$$
\bigg\|\sum^n_{j=1} c_j C_{\varphi_{j}}\bigg\|^2_e \geq \frac{1}{2\pi}\sum^n_{j=1} |c_j|^2 |J(\varphi_{i})|,
$$
where $|J(\varphi)|$ is the arclength measure of $J(\varphi)$. Accordingly, to study the questions discussed above, we assume throughout that our analytic self-maps $\varphi$ of $D$ satisfy $|\varphi(e^{i\theta})| < 1$ a.e.

\section{Preliminaries}
\setcounter{equation}{0}
\renewcommand{\theequation}{\thesection.\arabic{equation}}

Here we collect some preliminary facts used in the sequel.

\subsection{The Hardy and Bergman Spaces}

The Hardy space $H^2 = \calD_1$ is the set of all functions $f(z) = \sum^\infty_{n=0} a_nz^n$ analytic in $D$ with
$$
\|f\|^2 \equiv \sum^\infty_{n=0} |a_n|^2 < \infty.
$$
Given $f$ in $H^2$, the non-tangential limit $f(e^{i\theta}) = \lim_{\angle z \rarrow e^{i\theta}} f(z)$ exists for $d\theta$-almost every $e^{i\theta}$ in $\partial D$. Moreover, the correspondence $f(z) \rarrow f(e^{i\theta})$ allows one to think of $H^2$ as the closed subspace of $L^2 = L^2\left(\partial D, \frac{d\theta}{2\pi}\right)$ with orthonormal basis $\{e^{in \theta}\}^\infty_{n=0}$.

For $\alpha > -1$, the Bergman space $A^2_\alpha$ is the set of functions $f$ analytic in $D$ with
$$
\|f\|^2 = \frac{\alpha+1}{\pi} \int_D |f(z)|^2(1-|z|^2)^\alpha dA(z)
< \infty,
$$
where $dA$ is Lebesgue area measure on $D$. As mentioned above, $\calD_\beta = A^2_{\beta-2}$ for $\beta >1$, with equality of norms.

For information about $H^2$ and $A^2_\alpha$, see \cite{D} and \cite{CM}.

\subsection{Angular Derivatives} 

Let $\varphi$ be an analytic self-map of $D$. Then $\varphi$ has a (finite) angular derivative at $\zeta$ in $\partial D$ provided $\varphi(\zeta)$, the non-tangential limit of $\varphi$ at $\zeta$, exists and has modulus one, and
$$ \varphi'(\zeta) \equiv \lim_{\angle z \rarrow \zeta} \ \frac{\varphi(z)-\varphi(\zeta)}{z-\zeta}
$$
exists as a finite complex number. If the angular derivative $\varphi'(\zeta)$ fails to exist, we write $|\varphi'(\zeta)| = \infty$. In either case the Julia-Caratheodory Theorem \cite{CM} asserts in part that
$$
\liminf_{z \rarrow \zeta} \ \frac{1-|\varphi(z)|}{1-|z|} = |\varphi'(\zeta)|,
$$
where the limit inferior is taken unrestrictedly in $D$; moreover $|\varphi'(\zeta)| > 0$. Throughout we write $F(\varphi)$ for the set of all points in $\partial D$ where $\varphi$ has a finite angular derivative. For $\zeta$ in $F(\varphi)$ we have the relation $\varphi'(\zeta) = \overline{\zeta} \varphi(\zeta)|\varphi'(\zeta)|$. A condition necessary for the composition operator $C_\varphi$ to act compactly on $H^2$ is that $F(\varphi)$ be empty \cite{ST}. On the Bergman space $A^2_\alpha$, this condition is both necessary and sufficient \cite{MS}.

\subsection{Clark Measures}

Let $\varphi$ be an analytic self-map of $D$. If $|\alpha| =1$, there exists a finite positive Borel measure $\mu_\alpha$ on $\partial D$ such that
\begin{equation}
\frac{1-|\varphi(z)|^2}{|\alpha-\varphi(z)|^2} = \mbox{Re}\left(\frac{\alpha+ \varphi(z)}{\alpha-\varphi(z)}\right) = \int_{\partial D} P_z(e^{it})d\mu_\alpha (t)
\end{equation}
for $z$ in $D$, where
$$
P_z(e^{it}) = \frac{1-|z|^2}{|e^{it}-z|^2}
$$
is the Poisson kernel at $z$. The existence of $\mu_\alpha$ follows since the left side of equation (2.1) is a positive harmonic function. 
The measures $\mu_\alpha$ (the {\em Clark measures} of $\varphi$) were introduced as an operator-theoretic tool by D.~N.~Clark \cite{Cl}, and have been further analyzed by Alexsandrov \cite{A},
Poltoratski \cite{P}, and Sarason \cite{Sa2}.

On decomposing $\mu_\alpha = \mu_\alpha^{ac}+ \mu^s_\alpha$, where $\mu_\alpha^{ac}$ and $\mu_\alpha^s$ are, respectively, the absolutely continuous and singular parts with respect to Lebesgue measure, one finds by Fatou's theorem \cite{D} that
$$
\mu_{\alpha}^{ac} = \frac{1-|\varphi(e^{i\theta})|^2}{|\alpha-\varphi(e^{i\theta})|^2} \ \frac{d\theta}{2\pi}.
$$
The singular part $\mu^s_\alpha$ is carried by $\varphi^{-1}(\{\alpha\})$, the set of those $\zeta$ in $\partial D$ at which $\varphi(\zeta)$ exists and equals $\alpha$, and is itself the sum of the pure point measure
\begin{equation}
\mu_\alpha^{pp} = \sum_{\varphi(\zeta)=\alpha} \frac{1}{|\varphi'(\zeta)|} \ \delta_\zeta
\end{equation}
(here $\delta_\zeta$ is the unit point mass at $\zeta$) and a continuous singular measure $\mu_\alpha^{cs}$, either of which can vanish.

Let us write
$$
E(\varphi) = \overline{\dis\bigcup_{|\alpha|=1} spt(\mu^s_\alpha)},
$$
where $spt(\mu)$ denotes the closed support of a measure $\mu$. It is clear from Eqn.~(2.2) that $F(\varphi)$ is a subset of $E(\varphi)$.

\subsection{Essential Norms}

Let $H$ and $L$ be separable Hilbert spaces and write $\calB(H,L)$ for the space of bounded operators from $H$ to $L$. Let $\calK$ denote the subspace of compact operators in $\calB(H,L)$. The essential norm $\|T\|_e$ of an operator $T$ in $\calB(H,L)$ is the operator-norm distance from $T$ to $\calK$. We will find 
this alternate description useful:
\begin{equation}
\|T\|_e = \sup_{\{f_{n}\} \in \calU} \left(\limsup_{n \rarrow \infty} \|Tf_n\|\right),
\end{equation}
where $\calU$ is the collection of all sequences $\{f_n\}$ of unit vectors in $H$ which tend to zero weakly.

\subsection{Carleson Measures}

For a point $\zeta$ on the unit circle and $\delta > 0$, let
$S(\zeta,\delta)$ $= \{z \in D: \ |\zeta-z| < \delta\}$.
If $\mu$ is a finite positive Borel measure on $D$ and $\beta \geq 1$, we consider the quantities
\begin{equation}
\Delta_\beta(\mu)  =  \sup_{|\zeta|=1,\delta > 0}
\frac{\mu(S(\zeta,\delta))}{\delta^\beta}, \qquad
\Delta^*_\beta(\mu)  =  \limsup_{\delta \rarrow 0}
\left\{\sup_{|\zeta|=1} \frac{\mu(S(\zeta,\delta))}{\delta^\beta}\right\}.
\nonumber \end{equation}
One often says $\mu$ is a $\beta$-{\em Carleson} measure if $\Delta_\beta(\mu) < \infty$, and a {\em vanishing $\beta$-Carleson} measure if $\Delta^*_\beta(\mu) =0$.
The reader might consult \cite{CM} for the history of the following well-known result. The statement about $\|J\|_e$ can be deduced from ideas in the proof of Theorem 3.12 in 
\cite{CM}.

\renewcommand{\thetheorem}{\thesection.\arabic{theorem}}
\setcounter{theorem}{0}
\begin{theorem}
Let $\mu$ be a finite positive Borel measure on $D$ and assume $\beta \geq 1$. Then:

(i) The space $\calD_\beta$  (considered as a space of analytic functions on $D$) is contained in $L^2(\mu)$ if and only if $\mu$ is a $\beta$-Carleson measure. In this case the inclusion map $J: \calD_\beta \rarrow L^2(\mu)$ is bounded with norm comparable to $\sqrt{\Delta_\beta(\mu)}$.

(ii) If $\calD_\beta$ is contained in $L^2(\mu)$, then $\|J\|_e$, the essential norm of the inclusion map, is comparable to $\sqrt{\Delta_\beta^*(\mu)}$. In particular, $J$ is compact if and only if $\mu$ is a vanishing $\beta$-Carleson measure.
\end{theorem}

\subsection{A General Scheme for Compact Difference and Arc- \\ Connectedness}

For a bounded analytic function $w$ on $\partial D$, one can form the associated multiplication operator $M_w: f \rarrow wf$. If $\varphi$ is an analytic self-map of $D$, then we have the weighted composition operator $M_w C_\varphi$. Given two analytic self-maps of $D$, $\varphi$, and $\psi$, consider the self-maps $\varphi_t = t \varphi +(1-t)\psi$, $0 \leq t \leq 1$. Based on the formal operator identity
\begin{equation}
C_{\varphi_{s}}-C_{\varphi_{r}} = M_{\varphi-\psi}\left[\int^s_r C_{\varphi_{t}}dt\right]X,
\end{equation}
$0 \leq r < s \leq 1$, and the fact that the differentiation operator $X = \frac{d}{dz}$ is a topological isomorphism of $H^2_0$ (the subspace of $H^2$ consisting of all functions that vanish at the origin) and the Bergman space $A^2_1$ (see \cite{CM}), the second author and C.~Toews proved the following:

\begin{theorem}
{\rm \cite{MoT}} Let $\varphi$, $\psi$, and $\varphi_t$, $0 \leq t \leq 1$, be as above.

(i) Suppose the weighted composition operators $M_{\varphi-\psi}C_{\varphi_{t}}$ act boundedly from $A^2_1$ to $H^2$, with uniformly bounded norms, $0 \leq t \leq 1$. Then, there is a constant $B > 0$ such that, as operators on $H^2$,
$$
\|C_{\varphi_{s}}-C_{\varphi_{r}}\| \leq B|s-r|,
\ 0 \leq r < s < 1.
$$

(ii) Suppose that $\varphi$, $\psi$, and $\varphi_t$ satisfy the hypotheses of part (i) above, and in addition, that for each $t$, $0 \leq t \leq 1$, $M_{\varphi-\psi}C_{\varphi_{t}}$ is a compact operator from $A^2_1$ to $H^2$. Then $C_\varphi-C_\psi$ is a compact operator on $H^2$.
\end{theorem}

The above result remains true if one replaces $H^2$ and $A^2_1$ by $A^2_\alpha$ and $A^2_{\alpha+2}$, respectively, where $\alpha > -1$, see 
\cite{MoT}.

\section{Weighted composition operators on $H^2$ and $L^2$}
\setcounter{equation}{0}
\renewcommand{\theequation}{\thesection.\arabic{equation}}

For the analytic self-maps of $D$ considered here (those with $|\varphi(e^{i\theta})| < 1$ a.e.~on $\partial D$), Sarason 
\cite{Sa1} found a convenient representation of $C_\varphi$ as an integral operator on $H^2$ and even on the larger space $L^2 = L^2\left(\partial D, \frac{d\theta}{2\pi}\right)$. For $f$ in $L^2$, extend $f$ to a harmonic function in $D$ via the Poisson integral: $f(z) = \int_{\partial D} P_z(e^{it})f(e^{it}) \frac{dt}{2\pi}$. Putting $(C_\varphi f)(e^{i\theta}) = f(\varphi(e^{i\theta}))$ using the extended $f$ (since $|\varphi(e^{i\theta})| < 1$ a.e.), one has
\begin{equation}
(C_\varphi f)(e^{i\theta}) = \int_{\partial D} \frac{1-|\varphi(e^{i\theta})|^2}{|e^{it}-\varphi(e^{i\theta})|^2} \ f(e^{it}) \frac{dt}{2\pi}.
\end{equation}
Using the Schur test for boundedness of integral operators (stated below), Sarason showed that $C_\varphi$ is compact on $H^2$ if the Clark measures $\mu_\alpha$ of $\varphi$ are absolutely continuous for all $\alpha$ in $\partial D$; J.~H.~Shapiro and C.~Sundberg \cite{SSu2} established the converse via function-theoretic methods. Subsequently Cima and Matheson \cite{CiMa} discovered the following expression for the essential norm of an arbitrary $C_\varphi$ acting on $H^2$:
\begin{equation}
\|C_\varphi\|^2_e = \sup_{|\alpha| =1} \mu^s_\alpha(\partial D),
\end{equation}
a formula foreshadowed by C.~Cowen's inequalities for smooth $\varphi$ \cite[p.~84]{C}.

Here we adapt the integral operator approach  to investigate essential norms of {\em weighted} composition operators
$$
M_w C_\varphi: \ f \rarrow w \cdot(f\circ \varphi).
$$
We allow $w$ in $L^\infty$ and consider $M_w C_\varphi$ as mapping $L^2$ to $L^2$, $H^2$ to $L^2$ or (in the event that $w$ is in $H^\infty$), $H^2$ to $H^2$.
In all cases we have the following:

\begin{theorem}
Let $\varphi$ be an analytic self-map of $D$ with $|\varphi(e^{i\theta})| < 1$ a.e.~and 
having Clark measures $\mu_\alpha$, $|\alpha| = 1$. Suppose $w$ is a bounded measurable function on $\partial D$ such that $|w|$ is continuous at every point of $E(\varphi)$. Then
$$\dis\sup_{|\alpha|=1} \dis\int_{\partial D} |w|^2 d\mu^s_\alpha
\leq
\|M_w C_\varphi\|^2_e \leq 4
\dis\sup_{|\alpha|=1} \dis\int_{\partial D} |w|^2 d\mu^s_\alpha.
$$
In particular, $M_w C_\varphi$ is compact if and only if $w \equiv 0$ on $E(\varphi)$.
\end{theorem}

For the proof, our essential tool is the following.
\bigskip

\noindent{\bf The Schur Test} \cite[p.~282]{Nik} {\em Consider two measure spaces $(X,\mu)$ and $(Y,\nu)$, and let $N$ be a measurable function on the product space $Y \times X$. Suppose there exist positive measurable functions $p$ on $X$ and $q$ on $Y$ and constants $A, B > 0$ satisfying
$$
\int_X |N(y,x)|p(x)d\mu(x)  \leq  Aq(y), \ y \mbox{ in } Y,
\quad
\int_Y |N(y,x)|q(y)d\nu(y)  \leq  Bp(x), \ x \mbox{ in } X.
$$
Then the formula
$$
(Tf)(y) = \int_X N(y,x)f(x)d \mu (x),
$$
defines a bounded operator $T$ from $L^2(\mu)$ to $L^2(\nu)$ with $\|T\| \leq \sqrt{AB}$.}
\bigskip

We will also need several lemmas. The first, and the final conclusion in the second, are due to J.~E.~Shapiro \cite{Sh}.

\begin{lemma} {\rm \cite{Sh}}
Let $\varphi$ be an analytic self-map of $D$ with Clark measures $\mu_\alpha$, $|\alpha| =1$. If $f$ is continuous on $\partial D$, then
$$
\lim_{r \nearrow 1} \int_{\partial D}f \ \frac{1-r^2}{|\alpha-r\varphi|^2} \ \frac{d\theta}{2\pi} = \int_{\partial D} f \, d\mu^s_\alpha,
$$
for $|\alpha| =1$.
\end{lemma}

\begin{lemma}
If $f$ is continuous on $\partial D$, then
$$
\lim_{r \nearrow 1} \int_{\partial D}f \ \frac{1-|r\varphi|^2}{|\alpha-r\varphi|^2} \ \frac{d\theta}{2\pi} = \int_{\partial D} f \, d\mu_\alpha,
$$
uniformly in $\alpha$, $|\alpha| =1$. Thus $\int_{\partial D} f \, d\mu_\alpha$ is a continuous function
of $\alpha$.
\end{lemma}

{\em Proof.}
For $0 < r < 1$, the function $\frac{1-|r\varphi|^2}{|\alpha-r\varphi|^2}$ is bounded and harmonic on $D$, and thus is the Poisson integral of its boundary function. It follows from this, and the definition of $\mu_\alpha$, that the conclusion holds for $f = P_z$, the Poisson kernel at any $z$ in $D$. Thus the conclusion holds when $f$ is a finite linear combination of Poisson kernels. Such finite linear combinations are uniformly dense in the continuous functions on $\partial D$, and the lemma follows. \hfill $\Box$
\bigskip

The third lemma is a variant of exercise (7) in $\S$26 of Halmos' treatise \cite{H}; for the proof the interested reader can easily adapt the hint given there.

\begin{lemma}
Let $f$ and $f_n$, for $n = 1,2,3,\ldots$, be non-negative integrable functions on a measure space $(X,\mu)$ and suppose $\lambda \geq 0$. If $f_n \rarrow f$, $\mu$-a.e.~as $n \rarrow \infty$, and
$$
\limsup_{n \rarrow \infty}\left\{\int_X f_n d\mu - \int_X f \, d\mu\right\} \leq \lambda,
$$
then
$$
\limsup_{n \rarrow \infty}\int_X |f_n- f | d\mu \leq \lambda.
$$
\end{lemma}

{\em Proof of Theorem 3.1.} We can write $w = |w|v$ where $v$ is measurable and unimodular on $\partial D$. Since $M_w = M_vM_{|w|}$ and $M_v$ is unitary,
the theorem is unaltered by assuming that $w = |w| \geq 0$. For now we also assume that $w$ is continuous on $\partial D$.

Consider the normalized kernel function
$$
K_a(z) = \frac{\sqrt{1-|a|^2}}{1-\overline{a} z}, \qquad |a| < 1.
$$
If $|a| \rarrow 1$, then $K_a \rarrow 0$ weakly in $H^2$.
Now take $a = r \alpha$ where $0 \leq r < 1$ and $\alpha$ is fixed in $\partial D$. Then $|C_\varphi K_{r\alpha}|^2 = \frac{1-r^2}{|\alpha - r \varphi|^2}$ and
$$
\int |w|^2 d \mu^s_\alpha
 =  \lim_{r \nearrow 1} \int |w|^2 \frac{1-r^2}{|\alpha -r \varphi|^2} \ \frac{d\theta}{2\pi}
= \limsup_{r \nearrow 1} \|M_w C_\varphi K_{r\alpha}\|^2
$$
where Lemma 3.2 gives the first equality; the lower bound for $\|M_w C_\varphi\|^2_e$ follows from Eqn.~(2.3). (Throughout the proof, all integrals are taken over $\partial D$.)

For the upper bound we consider $M_w C_\varphi$ as an integral operator from $L^2 = L^2\left(\partial D, \frac{d\theta}{2\pi}\right)$ to $L^2\left(G, \frac{d\theta}{2\pi}\right)$, where $G = \{e^{i\theta}: w(e^{i\theta}) > 0\}$. By Eqn.~(3.1) the kernel of this operator is
$$
K(e^{i\theta},e^{it}) = w(e^{i\theta}) \, \frac{1-|\varphi(e^{i\theta})|^2}{|e^{it}-\varphi(e^{i\theta})|^2}.
$$
Similarly, if $0 < r < 1$, the kernel of the integral operator $M_w C_{r\varphi}$ is
$$
K_r(e^{i\theta},e^{it}) = w(e^{i\theta}) \, \frac{1-|r\varphi(e^{i\theta})|^2}{|e^{it}-r\varphi(e^{i\theta})|^2}.
$$
Since $\|r\varphi\|_\infty \leq r < 1$, $C_{r\varphi}$, and thus $M_w C_{r\varphi}$, are compact.
We apply the Schur test to the integral operator $M_w C_{r\varphi}-M_w C_\varphi$, which has integral kernel $N = K_r-K$. We take $\mu = \frac{d\theta}{2\pi}$, $\nu$ to be the restriction of $\frac{d\theta}{2\pi}$ to $G$, $p(e^{i\theta}) =1$ and $q(e^{i\theta}) = w(e^{i\theta})$.
Then
\begin{eqnarray*}
\lefteqn{\hspace*{-40pt}
\int |K_r(e^{i\theta},e^{it})-K(e^{i\theta},e^{it})|p(e^{it}) \frac{dt}{2\pi}
} \\[4mm]
& \leq & w(e^{i\theta}) \int
\left(\frac{1-|r\varphi(e^{i\theta})|^2}{|e^{it}-r\varphi(e^{i\theta})|^2}
+ \frac{1-|\varphi(e^{i\theta})|^2}{|e^{it}-\varphi(e^{i\theta})|^2}\right) \frac{dt}{2\pi}
\\[4mm]
& = & 2 w (e^{i\theta}) = 2q(e^{i\theta}),
\end{eqnarray*}
for all $e^{i\theta}$ in the circle. This is the first Schur hypothesis.

For the second Schur hypothesis, we write
$$
\lambda = \sup_{|\alpha| = 1} \int |w|^2 d\mu^s_\alpha,
$$
and consider sequences $r_n \nearrow 1$ and $\{\alpha_n\}$ in $\partial D$. It is enough to show that
\begin{equation}
\limsup_{n \rarrow \infty} \int |K_{r_{n}}(e^{i\theta},\alpha_n)-K(e^{i\theta},\alpha_n)|w(e^{i\theta}) \frac{d\theta}{2\pi} \leq 2 \lambda,
\end{equation}
for then, since $M_w C_{r_{n}\varphi}$ is compact, the Schur test with $A =2$ and $B=2\lambda$ will imply that
$$
\|M_wC_\varphi\|^2_e
 \leq  \limsup_{n \rarrow \infty} \|M_wC_{r_{n}\varphi}-M_wC_\varphi\|^2
 \leq  4\lambda,
$$
as desired.

We may assume that $\alpha_n$ tends to some $\alpha$ in $\partial D$ as $n \rarrow \infty$. First we use Lemma 3.4 with
$$
f_n = |w|^2 \frac{1-|\varphi|^2}{|\alpha_n-\varphi|^2}, \ f = |w|^2 \frac{1-|\varphi|^2}{|\alpha-\varphi|^2}.
$$
Note that
$$
\int f \frac{d\theta}{2\pi} = \int |w|^2 d\mu_\alpha - \int|w|^2 d\mu^s_\alpha,
$$
and similarly for $f_n, \mu_{\alpha_{n}}$ and $\mu^s_{\alpha_{n}}$, so that
$$
\int(f_n-f) \frac{d\theta}{2\pi} = \int|w|^2 d\mu_{\alpha_{n}} - \int|w|^2 d\mu_{\alpha} + \int |w|^2 d\mu^s_\alpha - \int|w|^2d\mu^s_{\alpha_{n}}.
$$
By the final conclusion in Lemma 3.3, the difference of the first two terms tends to zero as $n \rarrow \infty$, so that
$$
\limsup_{n \rarrow \infty} \int (f_n-f) \frac{d\theta}{2\pi} \leq \lambda.
$$
Since $f_n \rarrow f$ a.e., Lemma 3.4 implies that
\begin{equation}
\limsup_{n\rarrow \infty} \int |K(e^{i\theta}, \alpha_n)-K(e^{i\theta},\alpha)|w(e^{i\theta})\frac{d\theta}{2\pi} \leq \lambda.
\end{equation}
Now use Lemma 3.4 again, this time with
$$
f_n = |w|^2 \ \frac{1-|r_n\varphi|^2}{|\alpha_n-r_n \varphi|^2}, \ f = |w|^2 \ \frac{1-|\varphi|^2}{|\alpha-\varphi|^2}.
$$
We have
\begin{eqnarray*}
\int(f_n-f) \frac{d\theta}{2\pi} 
& \leq & \left|\int |w|^2 \
\frac{1-|r_n\varphi|^2}{|\alpha_n-r_n \varphi|^2}
\ \frac{d\theta}{2\pi} 
-\int |w|^2 d\mu_{\alpha_{n}}
\right|
\\[4mm]
& + & \left|
 \int |w|^2 d\mu_{\alpha_{n}}
- \int |w|^2 d\mu_{\alpha}
\right|
\\[4mm]
& + & \left|
\int |w|^2 d\mu_{\alpha}
- \int |w|^2 \
\frac{1-|\varphi|^2}{|\alpha-\varphi|^2} \
\frac{d\theta}{2\pi} 
\right|.
\end{eqnarray*}
As $n \rarrow \infty$, the first two terms on the right tend to zero by Lemma 3.3,
while the last term is exactly $\int|w|^2 d\mu^s_\alpha$, which does not exceed $\lambda$. Thus, according to Lemma 3.4,
$$
\limsup_{n \rarrow \infty} \int |K_{r_{n}}(e^{i\theta},\alpha_n) - K(e^{i\theta},\alpha)|w(e^{i\theta}) \frac{d\theta}{2\pi} \leq \lambda.
$$
Combining this with Eqn.~(3.4) yields Eqn.~(3.3) and thus the desired upper bound for $\|M_wC_\varphi\|^2_e$.

It remains to consider the case where $w = |w|$ is continuous at each point of $E(\varphi)$, but not necessarily on all of $\partial D$. Since $E(\varphi)$ is compact, there exists a function $g$ continuous on $\partial D$ with $g = w$ on $E(\varphi)$. Then $w-g$ is continuous at and vanishes at every point of $E(\varphi)$. Given $\epsilon >0$ one can construct a function $h$ continuous on $\partial D$ with $h = \epsilon$ on $E(\varphi)$ and $|w-g| \leq h$ on $\partial D$. Since
$\|M_{w-g}C_\varphi f\| \leq \|M_h C_\varphi f\|$ for all $f$ in $L^2$, $\|M_{w-g}C_\varphi\|_e \leq \|M_h C_\varphi\|_e$ by Eqn.~(2.3). An application of this inequality and the triangle inequality yields
$$
\|M_gC_\varphi\|_e - \|M_hC_\varphi\|_e \leq \|M_wC_\varphi\|_e 
\leq \|M_g C_\varphi\|_e + \|M_hC_\varphi\|_e.
$$
Because $g$ and $h$ are both continuous on $\partial D$, we can apply our earlier argument to estimate $\|M_gC_\varphi\|_e$ and $\|M_hC_\varphi\|_e$; in particular, $\|M_hC_\varphi\|_e = O(\epsilon)$ as $\epsilon \rarrow 0$. Since $\epsilon$ is arbitrary and $w=g$ on $E(\varphi)$, the theorem follows. 
\hfill $\Box$ \bigskip

\section{Local Compact Difference and a Sum Theorem for $H^2$}
\setcounter{equation}{0}
\renewcommand{\theequation}{\thesection.\arabic{equation}}

Our plan for studying linear combinations in $\calB(H^2)/\calK$ is to decompose, mod $\calK$, a composition operator into pieces associated to subsets of $E(\varphi)$ in a manner analogous to the decomposition in Theorem B above. This depends on an $H^2$ analogue of Theorem C above using Theorem 3.1 and $H^2$ versions of ideas from \cite{Mo}. Throughout, $\varphi$ and $\psi$ are analytic self-maps of $D$ with $|\varphi| < 1$ a.e. and $|\psi| <1$ a.e.~on $\partial D$, $\rho= \left|\frac{\varphi-\psi}{1-\overline{\varphi}\psi}\right|$, and $\varphi_t = t\varphi+(1-t)\psi$, $0 \leq t \leq 1$. We require three lemmas.

\begin{lemma}
Let $a < 1$ and suppose $G$ is a measurable subset of $\partial D$ with $\rho \leq a$ on $G$. Assume that $0 \leq t \leq 1$. Then $M_{\chi_{_{G}}}M_{\varphi-\psi}C_{\varphi_{t}}$ acts boundedly from $A^2_1$ to $L^2$, and
$$
\|M_{\chi_{_{G}}}M_{\varphi-\psi}C_{\varphi_{t}}\|_{\calB(A^2_{1},L^2)} \leq \frac{c}{1-a}
\|M_{\chi_{_{G}}}M_\rho C_{\varphi_{t}}\|_{\calB(H^2,L^2)},
$$
where $c$ is an absolute constant. Moreover, the same inequality  (with a different $c$) holds if both norms are replaced by the corresponding essential norms. In particular, if the operator $M_{\chi_{_{G}}}M_\rho C_{\varphi_{t}}: H^2 \rarrow L^2$ is compact, so is $M_{\chi_{_{G}}}M_{\varphi-\psi}C_{\varphi_{t}}: A^2_1 \rarrow L^2$.
\end{lemma}

{\em Proof.}
We consider the measures
$$
\nu_t = \left(\chi_{_{G}} |\varphi-\psi|^2\frac{d\theta}{2\pi}\right) \circ \varphi^{-1}_t, \ \beta_t = \left(\chi_{_{G}} \rho^2 \frac{d\theta}{2\pi}\right)\circ \varphi^{-1}_t
$$
on the disk $D$. As in \cite{Mo} we have
\begin{eqnarray*}
\frac{1-|\varphi_t|^2}{|1-\overline{\varphi} \psi|}
& = & \left|1+\overline{\varphi} \, 
\frac{(\psi-\varphi_t)}{1-\overline{\varphi} \psi} + \varphi_t
\frac{(\overline{\varphi}-\overline{\varphi}_t)}{1-\overline{\varphi} \psi}\right|
\\[4mm]
& = & \left|1-t \overline{\varphi} \,
\frac{(\varphi-\psi)}{1-\overline{\varphi} \psi} + (1-t)\varphi_t \,
\frac{(\overline{\varphi}-\overline{\psi})}{1-\overline{\varphi} \psi}\right|
\\[4mm]
& \geq & 1-\rho.
\end{eqnarray*}
Thus, if $\zeta$ lies in $\partial D$, $\delta > 0$
and $S(\zeta,\delta) = \{z \in D: |z-\zeta| < \delta\}$, on the set $G \cap \varphi^{-1}_t (S(\zeta,\delta))$ we have
$$
|\varphi-\psi|^2 = \rho^2|1-\overline{\varphi} \psi|^2 \leq \rho^2 \left(\frac{1-|\varphi_t|^2}{1-a}\right)^2 \leq \frac{4\delta^2}{(1-a)^2} \ \rho^2.
$$
It follows that
$$
\frac{\nu_t(S(\zeta,\delta))}{\delta^3} \leq \frac{4}{(1-a)^2} \ \frac{\beta_t (S(\zeta,\delta))}{\delta},
$$
so in the terminology of Section 2.5,
$$
\Delta_3(\nu_t) \leq \frac{4}{(1-a)^2} \ \Delta_1(\beta_t) \mbox{ and }
\Delta^*_3(\nu_t) \leq \frac{4}{(1-a)^2} \ \Delta^*_1(\beta_t).
$$
Since for all bounded analytic functions $h$ we have
$$
\int_D |h|^2d \nu_t = \|M_{\chi_{_{G}}} M_{\varphi-\psi} C_{\varphi_{t}} h\|^2_{L^2}
\mbox{ and } \int_D |h|^2 d \beta_t = \|M_{\chi_{_{G}}} M_\rho C_{\varphi_{t}} h\|^2_{L^2},
$$
the desired conclusion follows from Theorem 2.1 and the formula (2.3). \hfill $\Box$ \bigskip

\begin{lemma}
Let $G$ be an open set in $\partial D$ on which $\rho$ is essentially bounded away from one, and assume that $E(\varphi) \cap E(\psi) \cap G$ is closed. Then $E(\varphi_t) \cap G$ is contained in $E(\varphi) \cap E(\psi) \cap G$ for $0 < t < 1$.
\end{lemma}

{\em Proof.}
For $z$ in $D$ we have
$$
\frac{1-|\varphi_t(z)|}{1-|z|} \geq t \, \frac{1-|\varphi(z)|}{1-|z|} + (1-t)
\, \frac{1-|\psi(z)|}{1-|z|}.
$$
On letting $z \rarrow \zeta$ in $\partial D$, we see that if $\varphi_t$ has a finite angular derivative at $\zeta$, so do $\varphi$ and $\psi$, that is, $F(\varphi_t) \subset F(\varphi) \cap F(\psi)$. The opposite containment holds by linearity of the angular derivative and so $F(\varphi_t) = F(\varphi) \cap F(\psi).$

Suppose that $I$ is an open arc whose closure lies in $G$ and does not intersect both $E(\varphi)$ and $E(\psi)$. Let $\mu_{\alpha,t}$ be a Clark measure for $\varphi_t$. The point masses of $\mu_{\alpha,t}$, if any, are carried by $F(\varphi_t) \subset E(\varphi) \cap E(\psi)$, so that $\mu_{\alpha,t}$ puts no mass at the endpoints of $I$. A theorem of J.~E.~Shapiro \cite{Sh} then states that
\begin{equation}
\mu^s_{\alpha,t}(I) = \lim_{r \nearrow 1} \int_I \frac{1-r^2}{|\alpha-r\varphi_t|^2} \ \frac{d\theta}{2\pi},
\end{equation}
where $\mu^s_{\alpha,t}$ is the singular part of $\mu_{\alpha,t}$. It is enough to show that this quantity is zero, for then $spt\left(\mu^s_{\alpha,t}\right)$ cannot intersect $I$, and consequently neither can $E(\varphi_t)$.

For some $a < 1$ we have $\rho \leq a$ on $G$. Thus for $e^{i\theta}$ in $G$, $\psi(e^{i\theta})$ lies in the closed pseudo-hyperbolic disk with pseudo-hyperbolic radius $a$ and pseudo-hyperbolic center $\varphi(e^{i\theta})$. On noting the Euclidean center and radius of this (also Euclidean) disk \cite[p.~44]{CM}, one can verify that there is a positive constant $c$, depending only on $a < 1$, such that
$$
|\alpha-r \varphi_t(e^{i\theta})| \geq c \max\{|\alpha-r \varphi(e^{i\theta})|,|\alpha-r \psi(e^{i\theta})|\}
$$
for all $\alpha$ in $\partial D$, $e^{i\theta}$ in $G$ and $0 < r < 1$. 

Given any $\zeta$ in the closure $\overline{I}$ of the above arc $I$, there exists an open arc $A(\zeta)$ in $\partial\mathbb{D}$ containing $\zeta$ with either $\overline{A(\zeta)} \cap E(\varphi)$ or $\overline{A(\zeta)} \cap E(\psi)$ empty. The open cover $\{A(\zeta): \zeta \in \overline{I}\}$ for $\overline{I}$ has a finite subcover $\{A_1,\ldots,A_n,B_1,\ldots,B_m\}$ such that all of the sets $\overline{A}_i \cap E(\varphi)$, $i = 1,\ldots,n$ and $\overline{B}_j \cap E(\psi)$, $j = 1,\ldots,m$, are empty; of course, this subcover could consist only of $A_i$'s, or only of $B_j$'s. By Eqn.~(4.1) and the previous paragraph,
\begin{eqnarray*}
\mu^s_{\alpha,t}(I)
& = & \lim_{r \nearrow 1} \int_I \frac{1-r^2}{|\alpha-r \varphi_t|^2} \
\frac{d\theta}{2\pi}
\nonumber \\[3mm]
& \leq & \lim_{r \nearrow 1} \frac{1}{c} \left[
\sum^n_{i=1} \int_{A_{i}}  \frac{1-r^2}{|\alpha-r \varphi|^2} \
\frac{d\theta}{2\pi} + 
\sum^m_{j=1} \int_{B_{j}}  \frac{1-r^2}{|\alpha-r \psi|^2} \
\frac{d\theta}{2\pi}\right]
\nonumber \\[3mm]
& = & \frac{1}{c} \left[
\sum^n_{i=1} \mu^s_{\alpha,0}(A_i)+
\sum^m_{j=1} \mu^s_{\alpha,1}(B_j)\right]
\nonumber \\[2mm]
& = & 0,
\end{eqnarray*}
as desired.
\hfill $\Box$ \bigskip

Our third lemma localizes Theorem 2.2. The proof is as in \cite{MoT}, but now uses the operator Equation (2.5) left-multiplied by $M_{\chi_{_{G}}}$.

\begin{lemma}
Let $G$ be a measurable  subset of $\partial D$. Suppose that the weighted composition operators $M_{\chi_{_{G}}} M_{\varphi-\psi}C_{\varphi_{t}}$ act boundedly from $A^2_1$ to $L^2$ with norms uniformly bounded for $0 \leq t \leq 1$. If $M_{\chi_{_{G}}} M_{\varphi-\psi}C_{\varphi_{t}}$ is compact for $0 \leq t \leq 1$, then $M_{\chi_{_{G}}}(C_\varphi-C_\psi)$ is compact from $H^2$ to $L^2$.
\end{lemma}

We can now state our local compact difference theorem.

\begin{theorem}
Let $U$ be an open subset of $\partial D$ whose boundary intersects neither $E(\varphi)$ nor $E(\psi)$. Suppose that $\rho$ can be re-defined on a set of measure zero (if necessary), so that $\lim_{e^{i\theta} \rarrow \zeta} \rho(e^{i\theta}) = \rho(\zeta) = 0$ for every $\zeta$ in $E(\varphi) \cap E(\psi) \cap U$. Then $M_{\chi_{_{U}}}(C_\varphi-C_\psi)$ is a compact operator from $H^2$ to $L^2$.
\end{theorem}

{\em Proof.}
Since the intersection of $E(\varphi) \cap E(\psi)$ with $U$ must be compact, there is an open subset $G$ of $U$, containing this intersection and such that $\rho \leq \frac{1}{2}$ on $G$. We have
$$
M_{\chi_{_{U}}}(C_\varphi-C_\psi) = M_{\chi_{_{G}}}(C_\varphi-C_\psi) + M_{\chi_{_{U \setminus G}}} C_\varphi - M_{\chi_{_{U \setminus G}}} C_\psi.
$$
The last two operators on the right are compact by Theorem 3.1. According to Lemma 4.3, compactness of $M_{\chi_{_{G}}}(C_\varphi-C_\psi)$ will follow if we can show that the operators $M_{\chi_{_{G}}} M_{\varphi-\psi} C_{\varphi_{t}}$ map $A^2_1$ into $L^2$, are uniformly bounded in norm and are each compact. Lemma 4.1 gives uniform boundedness since
\begin{equation}
\|C_{\varphi_{t}}\|^2_{\calB(H^2)} \leq \frac{2}{1-|\varphi_t(0)|} \leq 2 \max 
\left\{\frac{1}{1-|\varphi(0)|},\frac{1}{1-|\psi(0)|}\right\},
\end{equation}
see \cite{CM}, while compactness follows from Lemmas 4.1, 4.2, and Theorem 3.1.
\hfill $\Box$ \bigskip

We close this section with a sum theorem for $H^2$. Here $\varphi,\varphi_1,\varphi_2,\ldots,\varphi_n$ are analytic self-maps of $D$ with $|\varphi| < 1$ a.e.~and $|\varphi_i|<1$ a.e.~on $\partial D$, $i = 1,2,\ldots,n$.

\begin{theorem}
Let $\varphi,\varphi_1,\ldots,\varphi_n$ be as above. Suppose that
\begin{itemize}
\item[(a)]
The sets $E(\varphi_i)$, $i = 1,\ldots,n$ are pairwise disjoint, and $E(\varphi)$ coincides
with $E(\varphi_1) \cup \cdots \cup E(\varphi_n)$, and
\item[(b)]
The functions $\rho_i = \left|\frac{\varphi-\varphi_i}{1-\overline{\varphi} \varphi_i}\right|$ can be altered on a set of measure zero in $\partial D$ (if necessary) to guarantee that $\lim_{e^{i\theta}\rarrow \zeta} \rho_i(e^{i\theta}) = \rho_i(\zeta) = 0$ for all $\zeta$ in $E(\varphi_i)$, $i = 1,\ldots,n$.
\end{itemize}
Then, as operators on $H^2$, $C_\varphi \equiv C_{\varphi_{1}} + \cdots + C_{\varphi_{n}} \  (\mbox{\rm mod } \calK).$
\end{theorem}

{\em Proof.}
Clearly, we can
find pairwise disjoint open sets $U_1, \ldots,U_n$ in $\partial D$ such that $U_i$ contains $E(\varphi_i)$ and $\rho_i \leq \frac{1}{2}$ on $U_i$, $i = 1,\ldots,n$. Let $G = \partial D \setminus \bigcup^n_{i=1} U_i$, so that $\partial D = G \cup U_1 \cup \cdots \cup U_n$. For a measurable subset $B$ of $\partial D$, let us write $M_B$ for $M_{\chi_{_{B}}}$. Then
\begin{eqnarray*}
\lefteqn{\hspace*{-30pt}
C_\varphi - C_{\varphi_{1}} - C_{\varphi_{2}} - \cdots -C_{\varphi_{n}} 
} \\[2mm]
& = & (M_G + M_{U_{1}} + \cdots + M_{U_{n}}) 
(C_\varphi - C_{\varphi_{1}} - \cdots -C_{\varphi_{n}})
\\[2mm]
& = & M_G C_\varphi - \sum^n_{i=1} M_G C_{\varphi_{i}}
\\[2mm]
& + & M_{U_{1}}(C_\varphi-C_{\varphi_{1}}) - \sum_{i\neq 1} M_{U_{1}} C_{\varphi_{i}}
\\[2mm]
& + & M_{U_{2}}(C_\varphi-C_{\varphi_{2}}) - \sum_{i\neq 2} M_{U_{2}} C_{\varphi_{i}}
\\[2mm]
& + & \vdots
\\[2mm]
& + & M_{U_{n}}(C_\varphi-C_{\varphi_{n}}) - \sum_{i\neq n} M_{U_{n}} C_{\varphi_{i}}.
\end{eqnarray*}
Theorem 3.1 tells us that every individual term on the right is compact except possibly $M_{U_{i}}(C_\varphi-C_{\varphi_{i}}), \ i = 1,\ldots,n$. These, however, must be compact by Theorem 4.4.
\hfill $\Box$ \bigskip

\section{Linear Combinations mod $\calK$}
\setcounter{equation}{0}

In this section we consider a class $\calS$ of analytic self-maps of $D$ for which definitive computations can be done. For $\varphi_1,\ldots,\varphi_n$ in $\calS$, we determine which linear combinations $c_1C_{\varphi_{1}} + \cdots + c_n C_{\varphi_{n}}$ are compact. We begin with some results which hold for arbitrary analytic self-maps of $D$.

\subsection{A First-Order Lower Bound for $\|c_1C_{\varphi_{1}} + \cdots+ c_n C_{\varphi_{n}}\|_e$}

Two analytic self-maps of $D$ have the {\em same first-order data} at $\zeta$ in $\partial D$ provided $\zeta$ lies in both $F(\varphi)$ and $F(\psi)$, $\varphi(\zeta) = \psi(\zeta)$ and $\varphi'(\zeta) = \psi'(\zeta)$. A special case
of a theorem of MacCluer \cite{M} states that if $C_{\varphi_{1}},\ldots,C_{\varphi_{n}}$ act on $\calD_\beta$, $\zeta$ is in $\partial D$, and no two of the maps $\varphi_1,\ldots,\varphi_n$ have the same first-order data at $\zeta$, then
$$
\|c_1 C_{\varphi_{1}} + \cdots + c_n C_{\varphi_{n}}\|^2_e \geq \sum^n_{k=1} |c_k|^2 \frac{1}{|\varphi'_k(\zeta)|^\beta}.
$$
Here it is understood that if $\zeta$ is not in $F(\varphi_k)$, then $|\varphi'_k(\zeta)| = \infty$. We will need a minor but useful extension of this result. For $\zeta$ in $\partial D$ and $M > 0$, let $\Gamma^1_{\zeta,M}$ denote the curve in $D$ given by $\frac{|z-\zeta|}{1-|z|^2} = M$, the boundary of a non-tangential approach region with vertex at $\zeta$. At the point $\zeta$ the sides of this region make angle $\theta$ with the radius to $\zeta$, where $2\cos \theta = \frac{1}{M}$. We will use the notation ``$\lim_{\Gamma^1_{\zeta,M}}$'' to indicate a limit taken as $z \rarrow \zeta$ along the 
starboard leg of $\Gamma^1_{\zeta,M}$ (that is, {\em counterclockwise}).

\begin{lemma}
Suppose $\varphi$ and $\psi$ are analytic self-maps of $D$. Then
$$
\begin{array}{l}
\dis\lim_{\Gamma^1_{\zeta,M}} \ \dis\frac{1-|z|^2}{1-\overline{\varphi(z)}\psi(z)}
\\[5mm]
\qquad 
=
\left\{\!\!\begin{array}{ll}
\dis\frac{2}{\overline{(1+i \tan \theta)}|\varphi'(\zeta)| + (1+ i \tan \theta)|\psi'(\zeta)|}
& \parbox[t]{106pt}{if $\zeta \in F(\varphi) \cap F(\psi)$ \\ and $\varphi(\zeta) = \psi(\zeta),$}
\\[7mm]
\, 0 & \mbox{otherwise.}
\end{array}\right.
\end{array}
$$
\end{lemma}

{\em Proof.}
The Schwarz inequality applied to the $H^2$ kernel functions gives
\begin{equation}
\frac{1-|z|^2}{|1-\overline{\varphi(z)} \psi(z)|} \leq 
\left(\frac{1-|z|^2}{1-|\varphi(z)|^2}\right)^{\frac{1}{2}}
\left(\frac{1-|z|^2}{1-|\psi(z)|^2}\right)^{\frac{1}{2}}.
\end{equation}
Both factors on the right are bounded in $D$, and if $\zeta$ is not in $F(\varphi)$, say, the first factor tends to zero as  $z \rarrow \zeta$. Thus we may assume that $\zeta$ lies in $F(\varphi) \cap F(\psi)$. It is also clear that if $\varphi(\zeta) \neq \psi(\zeta)$, the left side of (5.1) tends to zero as $z \rarrow \zeta$ along $\Gamma^1_{\zeta,M}$. 

In the remaining case, $\zeta$ lies in $F(\varphi) \cap F(\psi)$ and $\varphi(\zeta) = \psi(\zeta)$. For $z$ in $\Gamma^1_{\zeta,M}$, we have
$$
\frac{1-\overline{\varphi(z)}\psi(z)}{1-|z|^2}
= \frac{1-|\varphi(z)|^2}{1-|z|^2} + \overline{\varphi(z)} \ M \ \frac{\varphi(z)-\psi(z)}{|\zeta-z|}.
$$
Since $\Gamma^1_{\zeta,M}$ is nontangential at $\zeta$, the first term on the right tends to $|\varphi'(\zeta)|$ as $z \rarrow \zeta$ along $\Gamma^1_{\zeta,M}$. Moreover on the counterclockwise
leg of $\Gamma^1_{\zeta,M}$, $\zeta-z \sim \zeta e^{i\theta}|\zeta-z|$ as $z \rarrow \zeta$.  Recall that $\zeta \overline{\varphi(\zeta)} \varphi'(\zeta) = |\varphi'(\zeta)|$ and similarly for $\psi$. Since $Me^{i\theta} = \frac{1}{2} (1+i \tan \theta)$ and $\varphi(\zeta) = \psi(\zeta)$, we have
\begin{eqnarray*}
\lim_{\Gamma^1_{\zeta,M}} \ \frac{1-\overline{\varphi(z)} \psi(z)}{1-|z|^2}
& = & |\varphi'(\zeta)| + M \zeta \overline{\varphi(\zeta)}e^{i\theta}(\psi'(\zeta)-\varphi'(\zeta))
\\[2mm]
& = & \frac{1}{2} \left[(1+i \tan \theta)|\psi'(\zeta)| + \overline{(1+i \tan \theta)} |\varphi'(\zeta)|\right],
\end{eqnarray*}
as desired.
\hfill $\Box$ \bigskip

For $\zeta$ in $F(\varphi)$ we call the vector $D_1(\varphi,\zeta) = (\varphi(\zeta),\varphi'(\zeta))$ the first-order data of $\varphi$ at $\zeta$. Suppose we fix analytic self-maps of $D$, 
$\varphi_1,\ldots,\varphi_n$. For $\zeta$ in $\partial D$, we denote by $\calD_1(\zeta)$ the set of first-order data vectors at $\zeta$ associated to these self-maps:
$$
\calD_1(\zeta) = \{D_1(\varphi_j,\zeta): 1 \leq j \leq n \mbox{ and } \zeta \in F(\varphi_j)\}.
$$
Throughout Section 5 we write $F = F(\varphi_1) \cup \cdots \cup F(\varphi_n)$.

\begin{theorem}
With the above notation, let $C_{\varphi_{1}},\ldots,C_{\varphi_{n}}$ act on $\calD_\beta$. Then for any complex numbers $c_1,\ldots,c_n$ and $\zeta$ in $F$,
$$
\|c_1 C_{\varphi_{1}} + \cdots c_n C_{\varphi_{n}}\|^2_e \geq \sum_{\boldd \in \calD_{1}(\zeta)}
\bigg|\!\!\!\!\!\!
\sum_{\begin{array}{c}
\scriptstyle{\zeta \in F(\varphi_j)} \\[-1mm] \scriptstyle{D_1(\varphi_j,\zeta)=\boldd }
\end{array}}\!\!\!\!\!\! c_j\bigg|^2 \frac{1}{|d_1|^\beta},
$$
where $\boldd = (d_0,d_1)$.
\end{theorem}

{\em Proof.}
Referring to Lemma 5.1, we see that if $\varphi$ and $\psi$ are analytic self-maps of $D$, then $M$ tending to infinity means that $\theta \rarrow \frac{\pi}{2}$, so that
$$
\lim_{M \rarrow \infty}  \lim_{\Gamma^1_{\zeta,M}} 
\frac{1-|z|^2}{1-\overline{\varphi(z)}\psi(z)} = \left\{\!\!\!\begin{array}{ll}
\dis\frac{1}{|\varphi'(\zeta)|}
& \!\parbox[t]{120pt}{if $ \zeta \in F(\varphi) \cap F(\psi)$
and \\ $D_1(\varphi,\zeta) = D_1(\psi,\zeta),$} 
\\[7mm]
0 & \!\mbox{otherwise.}
\end{array}\!\!\!\right.
$$
As $|z| \rarrow 1$ the normalized kernel functions $\frac{k_z}{\|k_z\|}$ tend weakly to zero in $\calD_\beta$. Since $C^*_\varphi k_z = k_{\varphi(z)}$, we see from Eqn.~(2.3) that
\begin{eqnarray*}
\lefteqn{\hspace*{-30pt}
\|c_1 C_{\varphi_{1}} + \cdots + c_n C_{\varphi_{n}}\|^2_e
} \\[4mm]
& \geq & \lim_{M \rarrow \infty} \lim_{\Gamma^1_{\zeta,M}}
\left\|\left(\overline{c}_1 C^*_{\varphi_{1}} + \cdots + \overline{c}_n C^*_{\varphi_{n}}\right)
\frac{k_z}{\|k_z\|}\right\|^2_e
\\[4mm]
& = & \sum^n_{j,\ell=1} \overline{c}_j c_\ell \lim_{M \rarrow \infty} \lim_{\Gamma^1_{\zeta,M}}
\left(\frac{1-|z|^2}{1-\overline{\varphi_j(z)} \varphi_\ell(z)}\right)^\beta
\\[4mm]
& = & \!\!\! \sum_{\begin{array}{c}
\scriptstyle{\zeta \in F(\varphi_{j}) \cap F(\varphi_\ell)} 
\\
\scriptstyle{D_{1}(\varphi_j,\zeta) = D_{1}(\varphi_{\ell},\zeta)}
\end{array}}  \!\!\!\!
c_j \overline{c}_\ell \ \frac{1}{|\varphi_j(\zeta)|^\beta},
\end{eqnarray*}
which is a restatement of the desired conclusion. 
\hfill $\Box$ \bigskip

\begin{corollary}
If $c_1C_{\varphi_{1}} + \cdots + c_n C_{\varphi_{n}}$  is compact on $\calD_\beta$, then for every $\zeta$ in $F$ and every $\boldd$ in $\calD_1(\zeta)$,
$$
\sum_{\begin{array}{c}
\scriptstyle{\zeta \in F(\varphi_{j})} 
\\
\scriptstyle{D_{1}(\varphi_j,\zeta) = \boldd}
\end{array}}  \!\!\!\!
c_j =0.
$$
\end{corollary}

\subsection{A Remark on Theorem B}

Corollary 5.3 shows that the hypothesis on angular derivative sets in Theorem B is actually implied by condition (ii) of the theorem.

\begin{corollary}
Suppose that $C_\varphi, C_{\varphi_{1}}, \ldots, C_{\varphi_{n}}$ act on $\calD_\beta$ and
$$
C_\varphi \equiv C_{\varphi_{1}} +  \cdots + C_{\varphi_{n}} \ (\mbox{\rm mod }\calK).
$$
Then $F(\varphi_1),\ldots,F(\varphi_n)$ are pairwise disjoint and
$F(\varphi_1) \cup \cdots \cup F(\varphi_n)$ coincides with $F(\varphi)$.
\end{corollary}

{\em Proof.}
Let us write $\varphi = \varphi_0$, $c_0 = 1$, and $c_j = -1$ for $j = 1,\ldots,n$, so that
$$
c_0 C_{\varphi_{0}} + c_1 C_{\varphi_{1}} + \cdots + c_n C_{\varphi_{n}} \in \calK.
$$
By Corollary 5.3, each of the sets
$\{j: \ \zeta \in F(\varphi_j) \mbox{ and } D_1(\varphi_j,\zeta) = \boldd\}$
is either empty or contains exactly two elements, namely zero and an integer from $\{1,\ldots,n\}$. The conclusion follows. 
\hfill $\Box$ \bigskip

\subsection{Lower Bounds from Higher-Order Data}

We have considered the first-order data $D_1(\varphi,\zeta) = (\varphi(\zeta),\varphi'(\zeta))$ for $\zeta$ in $F(\varphi)$. In what follows, we look at higher-order data vectors
$$
D_k(\varphi,\zeta) = (\varphi(\zeta),\varphi'(\zeta),\varphi''(\zeta),\ldots,\varphi^{(k)}(\zeta))
$$
at points where the corresponding derivatives make sense. Specifically, we say $\varphi$ has $k^{th}$-order data at $\zeta$ in $\partial D$ if there exist complex numbers $b_0,b_1,\ldots,b_k$ with $|b_0|=1$ so that 
$$
\varphi(z) = b_0 +b_1(z-\zeta) + \cdots + b_k(z-\zeta)^k +o(|z-\zeta|^k)
$$ 
as $z \rarrow \zeta$ unrestrictedly in $D$. In this case $\lim_{\angle z \rarrow \zeta} \varphi^{(j)}(z)$ exists and equals $j!b_j$ for $j=1,\ldots,n$ (see, for example, the argument on p.~47 in \cite{Sa2}); we refer to this limit as $\varphi^{(j)}(\zeta)$. Since $|b_0| =1$, $\zeta$ is in $F(\varphi)$ and $b_1$ is the angular derivative $\varphi'(\zeta)$.

The model for this definition is of course a map which continues analytically across $\partial D$ near $\zeta$. Aside from the partial Taylor expansion we want our $\varphi$ to inherit another property of analyticity: order of contact. We say an analytic self-map $\varphi$ of $D$ has {\em order of contact} $c > 0$ at $\zeta$ if $|\varphi(\zeta)| =1$ and
$$
\frac{1-|\varphi(e^{i\theta})|^2}{|\varphi(\zeta)-\varphi(e^{i\theta})|^c}
$$
is essentially bounded above and away from zero as $e^{i\theta} \rarrow \zeta$. To clarify this and subsequent calculations, we map to the upper half-plane $\Omega = \{w: \mbox{ Im } w > 0\}$. For any $\alpha$ in $\partial D$, consider the conformal map $\tau_\alpha(z) = i \ \frac{\alpha -z}{\alpha+z}$, which takes $D$ onto $\Omega$ and $\alpha$ to 0. If $\varphi$ has finite angular derivative at $\zeta$, then $u = \tau_{\varphi(\zeta)} \circ \varphi \circ \tau_\zeta^{-1}$ is an analytic self-map of $\Omega$ having non-tangential limit $u(0) =0$ and finite angular derivative $u'(0) = |\varphi'(\zeta)|$. Suppose for the moment that $\varphi$ has an analytic continuation to a neighborhood of $\zeta$, so the same is true of $u$ at the origin. For $w$ near zero, $u(w) = \sum^\infty_{n=1} a_n w^n$ with $a_1 = u'(0)$. If we assume $|\varphi| < 1$ a.e.~on $\partial D$, then
$$
\mbox{Im } u(x) = \sum^\infty_{n=1} (\mbox{Im } a_n)x^n
$$
is positive for real $x$ near (but not equal to) zero. The smallest natural number $n$ with $a_n$ non-real must be even, say $n = 2m$, and so Im $u(x) \sim (\mbox{Im } a_{2m})x^{2m}$ as $x \rarrow 0$; moreover, Im $a_{2m} > 0$. If follows that near zero the image of $\mathbb{R}$ under $u$ is approximated by the curve $y = cx^{2m}$ for appropriate $c > 0$.
Further, if $\tau_\zeta(e^{i\theta}) = x$, then
$$
\mbox{Im } u(x) = \frac{1-|\varphi(e^{i\theta})|^2}{|\varphi(\zeta)+\varphi(e^{i\theta})|^2},
$$
and we find that $\frac{1-|\varphi(e^{i\theta})|^2}{|\varphi(\zeta)-\varphi(e^{i\theta})|^{2m}}$ tends to a positive number
as $e^{i\theta} \rarrow \zeta$, so that $\varphi$ has order of contact 
$2m$ at $\zeta$.
\bigskip

\addtocounter{theorem}{1}
{\bf Definition 5.5 \ }
We say an analytic self-map $\varphi$ of $D$ has {\em sufficient data} at $\zeta$ in $\partial D$ if
\begin{itemize}
\item[(i)]
$\varphi$ has finite angular derivative at $\zeta$;
\item[(ii)] $\varphi$ has order of contact $2m$ at $\zeta$ for some natural number $m$;\
\item[(iii)]
$\varphi$ has $2m^{th}$-order data at $\zeta$.
\end{itemize}
\bigskip

If $\varphi$ has sufficient data at $\zeta$ with order of contact $2m$ and $u = \tau_{\varphi(\zeta)} \circ \varphi \circ \tau_\zeta^{-1}$ is as above, then $u$ has an analogous expansion at the origin,
$$
u(w) = \sum^{2m}_{j=0} \frac{u^{(j)}(0)}{j!} \ w^j + o(|w|^{2m});
$$
here the derivatives can be realized as the non-tangential limits $\lim_{\angle w \rarrow 0} u^{(j)}(w)$, $j = 0,1,\ldots,2m$. Moreover, if $1 \leq k \leq 2m$, the $k^{th}$ order data vector $D_k(\varphi,\zeta)$ determines and is determined by the corresponding data $u'(0),\ldots,u^{(k)}(0)$ of $u$ at zero.

Given $\zeta$ in $\partial D$, a natural number $k \geq 2$ and $M > 0$, let $\Gamma^k_{\zeta,M}$ denote the locus of the equation $\frac{|\zeta-z|^k}{1-|z|^2} = M$ in $D$, a curve having ``order of contact $k$'' with $\partial D$ at $\zeta$. We write ``$\lim_{\Gamma^k_{\zeta,M}}$'' to indicate a limit taken as $z$ tends to $\zeta$ along $\Gamma^k_{\zeta,M}$.

\begin{lemma}
Suppose analytic self-maps of $D$, $\varphi$ and $\psi$, both have sufficient data at $\zeta$ in $\partial D$. Let $u = \tau_{\varphi(\zeta)} \circ \varphi \circ \tau_\zeta^{-1}$ and $v = \tau_{\psi(\zeta)} \circ \psi \circ \tau_\zeta^{-1}$, so that for $w$ near zero in $\Omega$,
$$
u(w)  =  \sum^{2m}_{j=1} \ \frac{u^{(j)}(0)}{j!} \ w^j + o(|w|^{2m}),
\qquad
v(w)  =  \sum^{2n}_{j=1} \ \frac{v^{(j)}(0)}{j!} \ w^j + o(|w|^{2n}),
$$
where $2m$ and $2n$ are the respective orders of contact of $\varphi$ and $\psi$ at $\zeta$. Then, if $k \geq 2$ and $M > 0$,
$$
\dis\lim_{\Gamma^k_{\zeta,M}} \
\dis\frac{1-|z|^2}
{1-\overline{\varphi(z)}\psi(z)} 
  = \left[ \: \overline{\left(\frac{u'(0)}{2} -i \frac{M}{2^{k-1}k!} u^{(k)}(0)\right)}
+ \left(\frac{v'(0)}{2} -i \frac{M}{2^{k-1}k!} v^{(k)}(0)\right) \, \right]^{-1}
$$
provided $k \leq 2m$, $k \leq 2n$, and $D_{k-1}(\varphi,\zeta) = D_{k-1}(\psi,\zeta)$, while the limit is zero otherwise.
\end{lemma}

{\em Proof.} 
First we assume that $2m \geq k$, $2n \geq k$ and $D_{k-1}(\varphi,\zeta) = D_{k-1}(\psi,\zeta)$. We put $\alpha = \varphi(\zeta) = \psi(\zeta)$ and observe from direct calculation that with $w = \tau_\zeta(z)$, 
$$
\frac{v(w)- \overline{u(w)}}{2i \mbox{ Im } w} = \frac{1-\overline{\varphi(z)} \psi(z)}{1-|z|^2}
\cdot \frac{|\zeta+z|^2}{(\overline{\alpha}+\overline{\varphi(z)})(\alpha+\psi(z))}.
$$
If we let $z \rarrow \zeta$ along the curve $\Gamma^k_{\zeta,M}$, $w = \tau_\zeta(z)$ tends to zero along its image in $\Omega$, which is a slight enough perturbation of the curve $\tilde{\Gamma}_k$ defined by the equation 
$\frac{|w|^k}{\mbox{Im } w} = \frac{M}{2^{k-2}}$
that the latter can be used to compute our limit.  That is,
$$
\lim_{\Gamma^k_{\zeta,M}} \ \frac{1-\overline{\varphi(z)} \psi(z)}{1-|z|^2} = 
\lim_{\tilde{\Gamma}^k} \
\frac{v(w)- \overline{u(w)}}{2i \mbox{ Im } w}.
$$
Moreover,
\begin{equation}
\frac{v(w)- \overline{u(w)}}{2i \mbox{ Im } w}
= \frac{\mbox{Im } u(w)}{\mbox{Im } w} +
\frac{v(w)- u(w)}{2i \mbox{ Im } w}.
\end{equation}
We write
$$
u(w) = \sum^{2m}_{j=1} a_j w^j + o(|w|^{2m}), \ v(w) 
= \sum^{2n}_{j=1} b_j w^j + o(|w|^{2n}),
$$
and put $w = re^{i\theta}$. Consider the first term on the right in Eqn.~(5.2).  Since $a_1,a_2,\ldots,a_{2m-1}$ are real and Im $a_{2m} > 0$, we have
\begin{eqnarray*}
\mbox{Im } u(w)
& = & a_1 r \sin \theta + a_2 r^2 \sin 2\theta + \cdots + a_{2m-1} r^{2m-1}
\sin(2m-1)\theta
\\[2mm]
& + & |a_{2m}| r^{2m} \sin (t_{2m} +2m\theta) + o(r^{2m}),
\end{eqnarray*}
where $a_{2m} = |a_{2m}|e^{it_{2m}}$ with $0 < t_{2m} < \pi$. We take $w$
in $\tilde{\Gamma}^k$ and divide by $\mbox{Im }w = r \sin \theta = \frac{2^{k-2}}{M} \ r^k$ to see that
\begin{equation}
\lim_{\tilde{\Gamma}_k} \ \frac{\mbox{Im } u(w)}{\mbox{Im } w} = a_1 + (\mbox{Im } a_k) \frac{M}{2^{k-2}}, \quad 2 \leq k \leq 2m,
\end{equation}
where we use the fact that $a_k$ is real if $k < 2m$.

Now consider the second term in Eqn.~(5.2). Since $D_{k-1}(\varphi,\zeta) = D_{k-1}(\psi,\zeta)$, we have $a_j =b_j$ for $j < k$. Thus for $w$ in $\tilde{\Gamma}_k$,
$$
\frac{v(w)-u(w)}{2i \mbox{ Im } w}
 =  -i \ \frac{M}{2^{k-1}} (b_k-a_k) e^{ik\theta} +o(1).
$$
Since $\theta \rarrow 0$ as $w \rarrow 0$ along $\tilde{\Gamma}_k$,
$$
\lim_{\tilde{\Gamma}_{k}} \
\frac{v(w)-u(w)}{2i \mbox{ Im } w}
= -i \ \frac{M}{2^{k-1}} (b_k-a_k).
$$
This equation, Eqn.~(5.3) and Eqn.~(5.2) give the desired result.

Suppose now that $2 \leq k \leq \min\{2m,2n\}$ and $D_{k-1}(\varphi,\zeta) \neq D_{k-1}(\psi,\zeta)$. Let $p$ be the smallest integer with $p \leq k-1$ and $a_p \neq b_p$. Then for $w = re^{i\theta}$ in $\tilde{\Gamma}_k$,
$$
\frac{v(w)-u(w)}{2i \mbox{ Im } w}
= -i \ \frac{M}{2^{k-1}} \left(\frac{b_p-a_p}{r^{k-p}}\right)e^{ip\theta} + o \left(\frac{1}{r^{k-p}}\right),
$$
a quantity whose modulus tends to infinity as $w \rarrow 0$ along $\tilde{\Gamma}_k$. Thus
\begin{equation}
\lim_{\Gamma^k_{\zeta,M}} \ \frac{1-|z|^2}{1-\overline{\varphi(z)} \psi(z)} = 0.
\end{equation}
Finally, if $k > 2m$, Eqn.~(5.4) follows from Eqn.~(5.1), the definition of ``order of contact'' and simple estimates applied to the Clark measure inequality
$$
\frac{1-|\varphi(z)|^2}{|\varphi(\zeta)-\varphi(z)|^2}
\geq \int_{\partial D} P_z \
\frac{1-|\varphi|^2}{|\varphi(\zeta)-\varphi|^2} \ \frac{d\theta}{2\pi}.
$$
The case $k > 2n$ is similar.
\hfill $\Box$ \bigskip

Now fix analytic self-maps $\varphi_1,\ldots,\varphi_n$ of $D$ and $\zeta$ in $F$. Assume that any $\varphi_j$ having finite angular derivative at $\zeta$ (that is, $\zeta$ is in $F(\varphi_j)$) in fact has sufficient data at $\zeta$. Given an integer $k \geq 2$, we write $\mathbb{M}_k(\zeta)$ for the set of those integers $j$, $1 \leq j \leq n$, for which $F(\varphi_j)$ contains $\zeta$ and the order of contact of $\varphi_j$ at $\zeta$ is at least $k$. Let us write $\calD_k(\zeta) = \{D_k(\varphi_j,\zeta): \ j \in \mathbb{M}_k(\zeta)\}$. 
We have a higher-order analogue of Theorem 5.2.

\begin{theorem}
Assume that $\varphi_1,\ldots,\varphi_n$ are analytic self-maps of $D$ as described above, and let $\zeta$ be in $F$. If $c_1,\ldots,c_n$ are complex, $k \geq 3$ and notation is as above,
$$
\|c_1 C_{\varphi_{1}} + \cdots +c_n C_{\varphi_{n}}\|^2_e \geq \sum_{\boldd \in \calD_{k-1}(\zeta)}
\bigg|\!\!\!\!\!\!
\sum_{\begin{array}{c}
\scriptstyle{j \in \mathbb{M}_{k}(\zeta)} \\[-1mm] \scriptstyle{D_{k-1}(\varphi_j,\zeta)=\boldd }
\end{array}}\!\!\!\!\!\! c_j\bigg|^2 \frac{1}{|d_1|^\beta},
$$
where $\boldd = (d_0,d_1,\ldots,d_{k-1})$ and each $C_{\varphi_{j}}$ acts on $\calD_\beta$.
\end{theorem}

{\em Proof.}
First consider analytic self-maps $\varphi$ and $\psi$ of $D$, each of which has sufficient data, with respective orders of contact $2m$ and $2n$, at a given $\zeta$ in $\partial D$. Let $u$ and $v$ be related to $\varphi$ and $\psi$ as in Lemma 5.6. Since $u^{(j)}(0)$ is real for $1 \leq j < 2m$ and Im $u^{(2m)}(0) > 0$ (and similarly for $v^{(j)}(0)$), we see from Lemma 5.6 that
$$
\lim_{M \rarrow \infty} \lim_{\Gamma^{k-1}_{\zeta,M}} \frac{1-|z|^2}{1-\overline{\varphi(z)} \psi(z)}
= \left\{\begin{array}{cl}
\dis\frac{1}{|\varphi'(\zeta)|} & \parbox[t]{130pt}{if $k \leq 2m, \ k \leq 2n,$ and \\ $D_{k-1}(\varphi,\zeta) = D_{k-1}(\psi,\zeta)$,}
\\
\\[1mm]
\ 0 & \mbox{otherwise.}
\end{array}\right.
$$
Proceeding as in the proof of Theorem 5.2, we find
\begin{eqnarray*}
\|c_1 C_{\varphi_{1}} + \cdots + c_n C_{\varphi_{n}} \|^2_e
& \geq & \sum^n_{j,\ell=1} \overline{c_j} c_\ell \lim_{M \rarrow \infty}
 \lim_{\Gamma^{k-1}_{\zeta,M}}
\left(\frac{1-|z|^2}{1-\overline{\varphi_j(z)} \varphi_\ell(z)}\right)^\beta
\\[4mm]
& = & \!\!\!\!\!\!
\sum_{\begin{array}{c}
\scriptstyle{j,\ell \in\mathbb{M}_{k}(\zeta)} \\[-1mm] \scriptstyle{D_{k-1}(\varphi_j,\zeta)=
D_{k-1}(\varphi_\ell,\zeta)}
\end{array}}
\!\!\!\!\!\!
\overline{c_j} c_\ell \ \frac{1}{|\varphi'_j(\zeta)|^\beta},
\end{eqnarray*}
which is the desired conclusion. 
\hfill $\Box$ \bigskip

The above theorem yields a higher-order version of MacCluer's lower bound for $\|C_\varphi-C_\psi\|_e$ in \cite{M}.

\begin{corollary}
Fix $\zeta$ in $\partial D$ and analytic self-maps $\varphi$ and $\psi$ of $D$,
both of which have sufficient data at $\zeta$, with respective orders of contact $2m$ and $2n$.
\begin{itemize}
\item[(i)]
If $n < m$, then $\|C_\varphi-C_\psi\|^2_e \geq \frac{1}{|\varphi'(\zeta)|^\beta}$.
\item[(ii)]
If $n = m$ and $D_{2m-1}(\varphi,\zeta) \neq D_{2m-1}(\psi,\zeta)$, then
$$
\|C_\varphi-C_\psi\|^2_e \geq \frac{1}{|\varphi'(\zeta)|^\beta} + 
\frac{1}{|\psi'(\zeta)|^\beta}.
$$
\end{itemize}
\end{corollary}

{\em Proof.}
Apply Theorem 5.7 with $\varphi_1 = \varphi$, $\varphi_2 = \psi$, $c_1 = 1$, and $c_2 = -1$.
\hfill $\Box$ \bigskip

We need a more delicate version of Theorem 5.7, which is conveniently expressed in terms of the following formalism.
For $\beta \geq 1$, $\calD^+_\beta$ will denote the reproducing kernel Hilbert space of functions on the right half-plane $\Omega^+ = \{z: \mbox{ Re } z > 0\}$ having the kernel functions $k_w^+(z) = (z+ \overline{w})^{-\beta}$, $w$ in $\Omega^+$. These spaces appear in the literature with various defining normalizations. When $\beta =1$, $\calD_\beta^+$ is the Hardy space on $\Omega^+$, see \cite{Ho}. For $\beta > 1$, $\calD^+_\beta$ is the weighted Bergman space of all functions $f$ analytic on $\Omega^+$ for which
$$
\|f\|^2_{\calD^+_{\beta}} \equiv \frac{2^{\beta-2}(\beta-1)}{\pi}
\int_{\Omega^+} |f(x+iy)|^2 x^{\beta-2} dx \, dy
$$
is finite, see \cite[p.~74]{Sai}. For us a key fact (used in Section  5.5) is this:
$\{k^+_w: \ w \in \Omega_+\}$ {\em is a linearly independent set in $\calD^+_\beta$.}

Now let $\mathbb{M}_k(\zeta)$ and $\calD_k(\zeta)$ be as defined prior to Theorem 5.7.

\begin{lemma}
Let $\varphi_1,\ldots,\varphi_n$ be analytic self-maps of $D$. Fix $\zeta$ in $F$ and suppose that if $F(\varphi_j)$ contains $\zeta$, then $\varphi_i$ has sufficient data at $\zeta$, $j = 1,\ldots,n$. Let $u_j$ be related to $\varphi_j$ as $u$ is related to $\varphi$ in Lemma 5.6. Then if $A > 0$ and $k$ is an even natural number,
\begin{eqnarray*}
\|c_1 C_{\varphi_{1}} + \cdots + c_n C_{\varphi_{n}} \|^2_e
& \geq & \sum_{\boldd \in \calD_{k-1}(\zeta)} 
\bigg\|\!\!\!\!\!\!
\sum_{\begin{array}{c}
\scriptstyle{j \in \mathbb{M}_{k}(\zeta)} \\[-1mm] \scriptstyle{D_{k-1}(\varphi_j,\zeta)=
\boldd}
\end{array}}\!\!\!\!\!\!
\overline{c_j} k^+_{\frac{u_j'(0)}{2}-i Au^{(k)}_{j}(0)}\bigg\|^2_{\calD^+_{\beta}}.
\end{eqnarray*}
\end{lemma}

{\em Proof.}
We use Lemma 5.6 with $M = A2^{k-1}k!$. Then we have
\begin{eqnarray*}
\lefteqn{
\|c_1 C_{\varphi_{1}} + \cdots + c_n C_{\varphi_{n}} \|^2_e
} \\[4mm]
& \geq & \lim_{\Gamma^k_{\zeta,M}} \bigg\|(\overline{c_1} C^*_{\varphi_{1}}+\cdots + \overline{c_n} C^*_{\varphi_{n}}) \frac{k_z}{\|k_z\|}\bigg\|^2_{\calD_{\beta}}
\\[4mm]
& \geq & \lim_{\Gamma^k_{\zeta,M}} \sum \overline{c_j}c_\ell
\left(\frac{1-|z|^2}{1-\overline{\varphi_j(z)} \varphi_\ell(z)}\right)^\beta
\\[4mm]
& = & 
\!\!\!\!\!
\sum_{
\begin{array}{c}
\scriptstyle{j,\ell \in \mathbb{M}_{k}(\zeta)} 
\\[-1mm] 
\scriptstyle{D_{k-1}(\varphi_j,\zeta)= D_{k-1}(\varphi_\ell,\zeta)}
\end{array}
}
\!\!\!\!\!\!\!\!\!\!\!
\overline{c_j}c_\ell 
\!
\left[ \
\overline{
\left(\frac{u_j'(0)}{2}  -iAu_j^{(k)}(0)\right)}
+
\left(
\frac{u_\ell'(0)}{2} - iAu_\ell^{(k)}(0)
\!\right)
\right]^{-\beta}
\\[4mm]
& = &
\sum_{\boldd \in \calD_{k-1}(\zeta)} 
\bigg\|
\!\!\!\!\!\!
\sum_{\begin{array}{c}
\scriptstyle{j \in \mathbb{M}_{k}(\zeta)} \\[-1mm] 
\scriptstyle{D_{k-1}(\varphi_j,\zeta)=\boldd}
\end{array}
}
\!\!\!\!\!\!
\overline{c_j} k^+_{\frac{u_j'(0)}{2}-i Au^{(k)}_{j}(0)}
\bigg\|^2_{\calD^+_{\beta}},
\end{eqnarray*}
as desired. Note that since $u'_j(0) > 0$, $A > 0$ and Im $u_j^{(k)}(0) \geq 0$, the complex number $\frac{u_j(0)}{2} -i Au_j^{(k)}(0)$ lies in $\Omega^+$.
\hfill $\Box$ \bigskip

\subsection{The Class $\calS$ and Making $\|C_\varphi-C_\psi\|_e$ Small}

For maps $\varphi$ and $\psi$ with sufficient data at a given $\zeta$ in $\partial D$, Corollary 5.8 describes two obstructions to $\|C_\varphi-C_\psi\|_e$ being small:
\begin{itemize}
\item[(a)]
unequal orders of contact at $\zeta$;
\item[(b)]
equal order of contact $2m$ but $D_{2m-1}(\varphi,\zeta) \neq D_{2m-1}(\psi,\zeta)$.
\end{itemize}
In this section we estimate $\|C_\varphi-C_\psi\|_e$ in the absence of these obstructions and characterize when it is zero. We work within the class $\calS$ of analytic self-maps $\varphi$ of $D$ for which $E(\varphi)$ is a finite set (so that $E(\varphi) = F(\varphi)$) and such that $\varphi$ has sufficient data at each point of $F(\varphi)$. For simplicity we restrict attention to composition operators on $\calD_1 = H^2$. We write $\Omega$ for the upper half-plane $\{z: \mbox{ Im } z >0\}$; the pseudo-hyperbolic metric $\Lambda$ on $\Omega$ is given by
$$
\Lambda(z,w) = \left|\frac{z-w}{z-\overline{w}}\right|.
$$
Note that $0 \leq \Lambda < 1$ on $\Omega \times \Omega$.
Recall from Section 5.3 that if $\varphi$ has sufficient data at $\zeta$ in $F(\varphi)$ with order of contact $2m$, and $u = \tau_{\varphi(\zeta)} \circ \varphi \circ \tau_\zeta^{-1}$, then $u^{(2m)}(0)$ lies in $\Omega$.

\begin{proposition}
Fix $\zeta$ in $\partial D$ and suppose that $\varphi$ and $\psi$ have sufficient data 
with respective orders of contact $2m$ and $2n$ at $\zeta$, and
moreover that $\varphi(\zeta) = \psi(\zeta)$.
Write $\rho = \left|\frac{\varphi-\psi}{1-\overline{\varphi}\psi}\right|$, $u = \tau_{\varphi(\zeta)} \circ \varphi \circ \tau_\zeta^{-1}$ and $v = \tau_{\psi(\zeta)} \circ \psi \circ \tau_\zeta^{-1}$.
\begin{itemize}
\item[(a)]
If $2m \neq 2n$, or $2m = 2n$ and $D_{2m-1}(\varphi,\zeta) \neq D_{2m-1}(\psi,\zeta)$, then $\rho(e^{i\theta}) \rarrow 1$ as $e^{i\theta} \rarrow \zeta$.
\item[(b)]
If $2m = 2n$ and $D_{2m-1}(\psi,\zeta) = D_{2m-1}(\varphi,\zeta)$, then 
$$
\lim_{e^{i\theta} \rarrow \zeta} \rho(e^{i\theta}) = \Lambda \left(u^{(2m)}(0),v^{(2m)}(0)\right)< 1.
$$
\item[(c)]
$\rho(e^{i\theta}) \rarrow 0$ as $e^{i\theta} \rarrow \zeta$ if and only if $2m = 2n$ and $D_{2m}(\varphi,\zeta) = D_{2m}(\psi,\zeta)$.
\end{itemize}
\end{proposition}

{\em Proof.}
First assume that the hypotheses of (b) hold.
Since $\psi(\zeta) = \varphi(\zeta)$, direct computation shows that if $\tau_\zeta(e^{i\theta}) =x$,
$$
\rho(e^{i\theta}) = \left|\frac{u(x)-v(x)}{u(x)-\overline{v(x)}}\right|.
$$
For $i = 1,2,\ldots,2m-1$, the MacLaurin coefficient $a_i$ for $u$ is real and equal to the corresponding coefficient $b_i$ for $v$. On the other hand, $a_{2m} = u^{(2m)}(0)/(2m)!$ and $b_{2m} = v^{(2m)}(0)/(2m)!$ lie in $\Omega$. Clearly,
$$
\lim_{x\rarrow 0} 
\left|\frac{u(x)-v(x)}{u(x)-\overline{v(x)}}\right| = \Lambda(a_{2m},b_{2m}),
$$
which is the desired conclusion for (b). 

Clearly (b) implies (c); the interested reader can easily verify (a).
\hfill $\Box$ 
\bigskip

For $\varphi$, $\psi$, and $\zeta$ as in Proposition 5.10, let us use the convention that $\rho(\zeta) = \lim_{e^{i\theta} \rarrow \zeta} \rho(e^{i\theta})$. The limit exists by Proposition 5.10 and for $\varphi$, $\psi$  in $\calS$ this convention, at worst, redefines $\rho$ on a finite set, leaving the multiplication operator $M_\rho$ unaltered. Note, however, that $\rho(\zeta)$ as just defined is in general not the same as the non-tangential limit $\lim_{\angle z \rarrow \zeta} \rho(z)$.

\begin{theorem}
Suppose that $\varphi$ and $\psi$ are in $\calS$ with $F(\varphi) = F(\psi) = F$, and that $C_\varphi$ and $C_\psi$ act on $H^2$. Suppose that at each $\zeta$ in $F$, $\varphi$ and $\psi$ have common order of contact $2m(\zeta)$ with $D_{2m(\zeta)-1}(\varphi,\zeta) = D_{2m(\zeta)-1}(\psi,\zeta)$ and moreover that $\rho(\zeta) \leq \frac{1}{2}$. 
Then
$$
\frac{1}{4} \max_{\zeta \in F} \frac{1}{|\varphi'(\zeta)|} \rho(\zeta)^2 \leq \|C_\varphi - C_\psi\|^2_e
\leq B \max_{\alpha \in \varphi(F)} \sum_{\varphi(\zeta)=\alpha} \frac{1}{|\varphi'(\zeta)|} \rho(\zeta)^2,
$$
where $B$ is an absolute constant and $\varphi(F) = \{\varphi(\zeta): \, \zeta \in F\}$. Moreover, the lower bound is valid without the assumption that $\rho(\zeta) \leq \frac{1}{2}$, $\zeta \in F$.
\end{theorem}

{\bf Note}: Suppose that $\varphi$ and $\psi$ satisfy the hypotheses of the theorem, and that $\mu_\alpha$ and $\nu_\alpha$ are their respective Clark measures. Then $\mu_\alpha^s$ and $\nu_\alpha^s$ are pure point, and since $\varphi$ and $\psi$ have the same first-order data,
$$
\nu_\alpha^s = \mu_\alpha^s = \!\!\sum_{\begin{array}{c}
\scriptstyle{\zeta \in F} \\[-1mm] \scriptstyle{\varphi(\zeta) = \alpha}
\end{array}} \!\!\frac{1}{|\varphi'(\zeta)|} \delta_\zeta.
$$
The upper bound for $\|C_\varphi - C_\psi\|^2_e$ in the theorem then takes the form
$$
B \sup_{|\alpha|=1} \int \rho^2 d\mu^s_\alpha.
$$
We expect that there is a lower bound of the same type.
\medskip

{\em Proof of Theorem 5.11.}
For $\zeta$ in $F$ write $u_\zeta = \tau_{\varphi(\zeta)} \circ \varphi \circ \tau_\zeta^{-1}$ and $v_\zeta = \tau_{\psi(\zeta)} \circ \psi \circ \tau_\zeta^{-1}$.
According to Lemma 5.9, if $\zeta$ is in $F$ and $A > 0$, 
$$
\|C_\varphi - C_\psi\|^2_e \geq \bigg\|k^+_{\frac{u'_{\zeta}(0)}{2}
 - i Au_{\zeta}^{(2m(\zeta))}(0)} -
k^+_{\frac{v'_{\zeta}(0)}{2} - i A v_{\zeta}^{(2m(\zeta))}(0)}\bigg\|^2_{H^2_{+}}.
$$
A calculation analogous to \cite[Lemma 9.11]{CM} shows that for $z$, $w$ in $\Omega_+$,
$$
\|k_z^+-k_w^+\|^2_{H^2_{+}} = \left|\frac{z-w}{z+\overline{w}}\right|^2
\left(\frac{1}{z+\overline{z}} + \frac{1}{w+ \overline{w}}\right).
$$
We put
$$
A = \frac{u'_\zeta(0)}{\left|u_\zeta^{(2m(\zeta))}(0)-
\overline{v_\zeta^{(2m(\zeta))}(0)}\right|},
$$

$$ z = \frac{u'_\zeta(0)}{2} -iAu_\zeta^{(2m(\zeta))}(0), \ w = \frac{v'_\zeta(0)}{2} -iA v_\zeta^{(2m(\zeta))}(0),
$$
and note that $|z + \overline{w}| \leq 2u'_\zeta(0)$ while $|z-w| = u_\zeta'(0) \Lambda \left(u_\zeta^{(2m(\zeta))}(0),v_\zeta^{(2m(\zeta))}(0)\right)$. It follows that
$$
\|k_z^+-k_w^+\|^2_{H_{+}^2} \geq \frac{1}{4u_\zeta'(0)} \
\Lambda \left(u_\zeta^{(2m(\zeta))}(0),v_\zeta^{(2m(\zeta))}(0)\right)^2,
$$
which, by Proposition 5.10 and since $u_\zeta'(0) = |\varphi'(\zeta)|$ gives the lower bound for 
the essential norm of $C_\varphi-C_\psi$.

For the upper bound select a finite union $G$ of pairwise disjoint open arcs which
contains $F$ and is such that $\rho \leq \frac{2}{3}$ on $G$. Now
$$
C_\varphi-C_\psi = M_{\chi_{_{G}}}(C_\varphi-C_\psi)+M_{\chi_{_{\partial D\setminus G}}}C_\varphi-M_{\chi_{_{\partial D\setminus G}}}C_\psi,
$$
 and the last two terms on the right are compact operators by Theorem 3.1.
 Therefore $\|C_\varphi-C_\psi\|_e = \|M_{\chi_{_{G}}}(C_\varphi-C_\psi)\|_e$ which we now estimate. To use Eqn.~(2.3) it is permissible to restrict the sequence $\{f_n\}$ to a subspace of finite codimension. With notation as in Section 2.6, we pick a sequence of unit vectors $f_n$ in $H^2_0 = (\mbox{ker }X)^\perp$ which converges weakly to zero. When applied to an $H^2$ function, the operator identity (2.5) holds pointwise in $D$ and thus, since the kernel functions span $H^2$, in the weak operator topology. Thus
$$
\|M_{\chi_{_{G}}}(C_\varphi-C_\psi)f_n\|_{L^2}
\leq \|X\|_{\calB(H^2_0,A^2_{1})} \int^1_0 \|M_{\chi_{_{G}}} M_{\varphi-\psi}C_{\varphi_{t}}g_n\|_{L^2}dt,
$$
where $g_n = Xf_n/\|Xf_n\|$. Since $X$ is bounded below on $H^2_0$, $g_n \rarrow 0$ weakly in $A^2_1$. Now $\|M_{\chi_{_{G}}}M_{\varphi-\psi}C_{\varphi_{t}}\|_{\calB(A_{1}^2,L^2)}$ is bounded for $0 \leq t \leq 1$ by Lemma 4.1 and Eqn.~(4.2), so we can apply Fatou's Lemma to conclude that
$$
\limsup_{n\rarrow \infty} \|M_{\chi_{_{G}}}(C_\varphi-C_\psi)f_n\|_{L^2}
\leq \|X\|_{B(H^2_{0},A^2_{1})} \int^1_0 \limsup_{n\rarrow \infty} \|M_{\chi_{_{G}}}M_{\varphi-\psi}C_{\varphi_{t}}g_n\|_{L^2}dt.
$$
Moreover, by Eqn.~(2.3) and Lemma 4.1 we have
\begin{eqnarray*}
\limsup_{n\rarrow \infty} \|M_{\chi_{_{G}}}M_{\varphi-\psi}C_{\varphi_{t}}g_n\|_{L^2}
& \leq & \|M_{\chi_{_{G}}}M_{\varphi-\psi}C_{\varphi_{t}}\|_{e,\calB(A_{1}^2,L^2)}
\\[2mm]
& \leq & b\|M_{\chi_{_{G}}}M_{\rho}C_{\varphi_{t}}\|_{e,\calB(H^2,L^2)}
\\[2mm]
& \leq & b\|M_{h}C_{\varphi_{t}}\|_{e,\calB(H^2,L^2)},
\end{eqnarray*}
where $b>0$, $h$ is any continuous function on $\partial D$ with $h \geq \chi_{_{G}} \rho$ on $\partial D$ and $h(\zeta) = \rho(\zeta)$ for $\zeta$ in $F$. Clearly we can choose such an $h$ while also requiring that $\rho \leq \frac{3}{4}$ on the set $W = \{e^{i\theta}: \, h(e^{i\theta}) > 0\}$. Then, by Theorem 3.1 and Lemma 3.2, if $\mu_{\alpha,t}$ is a Clark measure for $\varphi_t$,
$$
\|M_h C_{\varphi_{t}}\|_e
 \leq  2 \sup_{|\alpha|=1} \left\{\int_W h^2 d\mu_{\alpha,t}^s\right\}^{\frac{1}{2}}
= 2\sup_{|\alpha|=1} \left\{\lim_{r \nearrow 1} \int_W h^2 
\frac{1-r^2}{|\alpha -r \varphi_t|^2} \ \frac{d\theta}{2\pi}\right\}^{\frac{1}{2}}.
$$
Since $\rho \leq \frac{3}{4}$ on $W$ we see (as noted in the proof of Lemma 4.2) that there is a constant $c > 0$ such that $|\alpha -r\varphi_t(e^{i\theta})| \geq c |\alpha-r \varphi(e^{i\theta})|$,  for $\alpha$ in $\partial D$, $e^{i\theta}$ in $W$, $0 < r < 1$ and $0 \leq t < 1$. Thus
\begin{eqnarray*}
\|M_h C_{\varphi_{t}}\|_e
& \leq & \frac{2}{c} \sup_{|\alpha|=1} \left\{\lim_{r \nearrow 1} \int_{\partial D} h^2 
\frac{1-r^2}{|\alpha-r\varphi|^2}  \frac{d\theta}{2\pi}\right\}^{\frac{1}{2}}
\\[4mm]
& = & \frac{2}{c} \sup_{|\alpha|=1} \int_{\partial D} h^2 
d\mu^s_\alpha ,
\end{eqnarray*}
where $\{\mu_\alpha\}$  are the Clark measures for $\varphi$. Since 
each $\mu^s_\alpha$ is pure point in the present circumstances,
$$
\mu^s_\alpha = \left\{\!\!\!\!\begin{array}{ll}
\dis\sum_{
\begin{array}{c}
\scriptstyle{\zeta \in F} \\[-1mm]
\scriptstyle{\varphi(\zeta) = \alpha }
\end{array}
}
\!\!\! \dis\frac{1}{|\varphi'(\zeta)|} \, \delta_\zeta
& \mbox{if } \alpha \in \varphi(F),
\\
\\
\quad 0 & \mbox{otherwise.}
\end{array}\right.
$$
In view of Eqn.~(2.3) and the arbitrariness of $\{f_n\}$ in $H^2_0$, this gives the desired upper bound for $\|C_\varphi-C_\psi\|_e$.
\hfill $\Box$ \bigskip

\begin{corollary}
Suppose that $\varphi$ and $\psi$ lie in $\calS$ with respective orders of contact $2m(\zeta)$ and $2n(\zeta)$ at each $\zeta$ in (respectively) $F(\varphi)$ and $F(\psi)$. Then the following are equivalent.
\begin{itemize}
\item[(i)]
$C_\varphi-C_\psi$ is compact on $H^2$.
\item[(ii)]
$F(\psi) = F(\varphi)$ (we call this set $F$) and for all $\zeta$ in $F$,  $2n(\zeta) = 2m(\zeta)$ and $D_{2m(\zeta)}(\psi,\zeta) = D_{2m(\zeta)}(\varphi,\zeta)$.
\end{itemize}
\end{corollary}

\subsection{Linear Relations Mod $\calK$}

We fix $\varphi_1,\ldots,\varphi_n$ in $\calS$ and 
again write $F$ for the union $F(\varphi_1) \cup \cdots \cup F(\varphi_n)$, a finite set. For $\zeta$ in $F$ and $k = 2,4,6,\ldots,$ let 
$$
\mathbb{N}_k(\zeta) = \{j: \ F(\varphi_j) \mbox{ contains $\zeta$ and $k$ is 
the order of contact of $\varphi_j$ at $\zeta\}$.}
$$
We also write $\calE_k(\zeta) = \{D_k(\varphi_j,\zeta): \ j \mbox{ is in } \mathbb{N}_k(\zeta)\}$.

\begin{theorem}
Let $\varphi_1,\ldots,\varphi_n$ be in $\calS$ and set notation as above. Given complex numbers $c_1,\ldots,c_n$, the following are equivalent:
\begin{itemize}
\item[(i)]
$c_1 C_{\varphi_{1}} + \cdots + c_n C_{\varphi_{n}}$ is compact on $\calD_\beta$;
\item[(ii)]
for every $\zeta$ in $F$, every even $k \geq 2$ and every $\boldd$ in $\calE_k(\zeta)$,
$$
\sum_{\begin{array}{c}
\scriptstyle{j \in \mathbb{N}_{k}(\zeta)}
\\[-1mm]
\scriptstyle{D_k(\varphi_j,\zeta) =\boldd}
\end{array}}
\!\!\! c_j = 0.
$$
\end{itemize}
\end{theorem}

{\em Proof.}
First assume that $c_1 C_{\varphi_{1}} + \cdots + c_n C_{\varphi_{n}}$ is compact. Fix $\zeta$ in $F$ and let $u_j$ be related to $\varphi_j$ as $u$ is related to $\varphi$ in Lemma 5.6. Let $\mathbb{M}_k(\zeta)$ and $\calD_k(\zeta)$ be as defined prior to Theorem 5.7. According to Lemma 5.9 (with $A=1$), if $k \geq 2$, $\boldd$ is in $\calD_{k-1}(\zeta)$, and $w_j = \frac{u_j'(0)}{2} - i u_j^{(k)}(0)$, then
$$
\sum_{\begin{array}{c}
\scriptstyle{j \in \mathbb{M}_{k}(\zeta)}
\\
\scriptstyle{D_{k-1}(\varphi_j,\zeta) =\boldd}
\end{array}}
\!\!\! \overline{c_j} k^+_{w_{j}} = 0.
$$
But if $D_{k-1}(\varphi_j,\zeta)$ and $D_{k-1}(\varphi_\ell,\zeta)$ coincide, $w_j = w_\ell$ exactly when $D_k(\varphi_j,\zeta) = D_k(\varphi_\ell,\zeta)$. Using the fact that $\{k^+_w: \mbox{ Re } w > 0\}$ is linearly independent
in $\calD^+_\beta$, we see that for every $\boldd$ in $\calD_k(\zeta)$,
\begin{equation}
\sum_{\begin{array}{c}
\scriptstyle{j \in \mathbb{M}_{k}(\zeta)}
\\
\scriptstyle{D_k(\varphi_j,\zeta) =\boldd}
\end{array}}
\!\!\! c_j = 0.
\end{equation}
Since $\mathbb{M}_k(\zeta)$ is the union of the disjoint sets $\mathbb{N}_k(\zeta)$ and $\mathbb{M}_{k+1}(\zeta)$, we see from Eqn.~(5.5) that
$$
\sum_{\begin{array}{c}
\scriptstyle{j \in \mathbb{N}_{k}(\zeta)}
\\
\scriptstyle{D_k(\varphi_j,\zeta) =\boldd}
\end{array}}
\!\!\! c_j + \!\!\!\!
\sum_{\begin{array}{c}
\scriptstyle{j \in \mathbb{M}_{k+1}(\zeta)}
\\
\scriptstyle{D_k(\varphi_j,\zeta) =\boldd}
\end{array}}
\!\!\! c_j = 0.
$$
Suppose now $\boldd = (d_0,d_1,\ldots,d_k)$. There are a finite number of elements in $\calD_{k+1}(\zeta)$ having the form $(d_0,d_1,\ldots,d_k,*)$, call them $\boldd_i = (d_0,d_1,\ldots,d_k,d^i_{k+1})$, $i = 1,2,\ldots,r$. Then
$$
\sum_{\begin{array}{c}
\scriptstyle{j \in \mathbb{M}_{k+1}(\zeta)}
\\
\scriptstyle{D_k(\varphi_j,\zeta) =\boldd}
\end{array}}
\!\!\! c_j = \sum^r_{i=1} \bigg\{\!\!\!\!\!
\sum_{\begin{array}{c}
\scriptstyle{j \in \mathbb{M}_{k+1}(\zeta)}
\\
\scriptstyle{D_{k+1}(\varphi_j,\zeta) =\boldd_i}
\end{array}}
\!\!\!\!\!\! c_j \bigg\}
$$
and each of the $r$ summands on the right vanishes, again by Eqn.~(5.5). Therefore (ii) holds.

Now assume that (ii) holds. We consider the case $\beta =1$. The set $F$ is finite, say $F = \{\zeta_1,\ldots,\zeta_s\}$. Let $U_1,\ldots,U_s$ be disjoint arcs whose union is $\partial D$ and for which the interior of $U_i$ contains $\zeta_i$ for $i = 1,\ldots,s$. Clearly,
$$
\sum^n_{j=1} c_j C_{\varphi_{j}} = \sum^s_{i=1} M_{\chi_{U_{i}}} \left(\sum^n_{j=1} c_j C_{\varphi_{j}}\right).
$$
Taking $U$ to be any $U_i$ and $\zeta = \zeta_i$, it is enough to show that the operator
$$
M_{\chi_{_{U}}}\left(\sum^n_{j=1} c_j C_{\varphi_{j}}\right),
$$
considered as acting from $H^2$ to $L^2$, is compact. This operator can be written as
\begin{equation}
\sum_{\zeta \notin F(\varphi_j)} c_j M_{\chi_{_{U}}} C_{\varphi_{j}} + \sum_{\zeta \in F(\varphi_{j})} c_j 
M_{\chi_{_{U}}} C_{\varphi_{j}}.
\end{equation}
If $F(\varphi_j)$ does not contain $\zeta$, then $\chi_{_{U}}$ is identically zero in a neighborhood of $F(\varphi_j)$, and we see from Theorem 3.1 that $M_{\chi_{_{U}}}C_{\varphi_{j}}$ is compact. Thus the first sum in the expression (5.6) is a compact operator.

Consider now the second sum, which can be rewritten as
$$
\sum_{m \geq 1} \left\{\sum_{j \in \mathbb{N}_{2m}(\zeta)}
c_j M_{\chi_{_{U}}} C_{\varphi_{j}}\right\};
$$
the sum over $m$ is of course finite since $\mathbb{N}_{2m}(\zeta)$ is empty for $m$ large enough. For nonempty $\mathbb{N}_{2m}(\zeta)$, 
$$
\sum_{j \in  \mathbb{N}_{2m}(\zeta)} c_j M_{\chi_{_{U}}}C_{\varphi_{j}}
= \sum_{\boldd \in \calE_{2m}(\zeta)}
\Bigg\{\!\!\!\!\!\sum_{\begin{array}{c}
\scriptstyle{j \in \mathbb{N}_{2m}(\zeta)}
\\[-1mm]
\scriptstyle{D_{2m}(\varphi_j,\zeta) =\boldd}
\end{array}} \!\!\!
c_j M_{\chi_{_{U}}} C_{\varphi_{j}}\Bigg\}.
$$
If $\boldd$ is in $\calE_{2m}(\zeta)$, $j$ and $\ell$ are in $\mathbb{N}_{2m}(\zeta)$ and $D_{2m}(\varphi_j,\zeta) = D_{2m}(\varphi_\ell,\zeta)= \boldd$, Proposition 5.10(c)  tells us that
$$
\left|\frac{\varphi_j(e^{i\theta})-\varphi_\ell(e^{i\theta})}{1-\overline{\varphi_j(e^{i\theta})}
\varphi_\ell(e^{i\theta})}\right| \rarrow 0 \mbox{ as } e^{i\theta} \rarrow \zeta.
$$
It follows from Theorem 4.4 that $M_{\chi_{_{U}}}(C_{\varphi_{j}}-C_{\varphi_{\ell}})$ is compact. Fix such $\ell$ and write $T = M_{\chi_{_{U}}}C_{\varphi_{\ell}}$.  Then for each $j$ in $\mathbb{N}_{2m}(\zeta)$ with $D_{2m}(\varphi_j,\zeta) = \boldd$, there is a compact operator $K_j$ with $M_{\chi_{_{U}}}C_{\varphi_{j}} = T+K_j$. Thus
$$
\!\!\!
\sum_{\begin{array}{c}
\scriptstyle{j \in \mathbb{N}_{2m}(\zeta)}
\\[-1mm]
\scriptstyle{D_{2m}(\varphi_j,\zeta) =\boldd}
\end{array}}
\!\!\!
c_j M_{\chi_{_{U}}} C_{\varphi_{j}}
= \bigg(\!\!\!\!\sum_{\begin{array}{c}
\scriptstyle{j \in \mathbb{N}_{2m}(\zeta)}
\\[-1mm]
\scriptstyle{D_{2m}(\varphi_j,\zeta) =\boldd}
\end{array}}\!\!\! c_j
\bigg)T +  \!\!\!
\sum_{\begin{array}{c}
\scriptstyle{j \in \mathbb{N}_{2m}(\zeta)}
\\[-1mm]
\scriptstyle{D_{2m}(\varphi_j,\zeta) =\boldd}
\end{array}}
\!\!\! c_j K_j.
$$
The coefficient of $T$ vanishes by the hypothesis (ii), and we are left with a compact operator,
verifying (i). The proof for the case $\beta > 1$ is similar, with Theorem C playing the role of Theorem 4.4, the sets $U_i$ taken to lie in $D$, and Proposition 5.10(c) replaced by the assertion that if $2m=2n$ and $D_{2m}(\varphi,\zeta) = D_{2m}(\psi,\zeta)$, then $\rho(z) \rarrow 0$ as $z \rarrow \zeta$ unrestrictedly in $D$, an implication easily established by calculations in the proofs of Lemma 5.6 and Proposition 5.10.  
\hfill $\Box$ 

\begin{remark} {\rm
It is sometimes convenient to rephrase condition (ii) in Theorem 5.13. With notation as in the statement, fix $\zeta$ in $F$, an even natural number $k$ and a vector $\boldd$ in $\calE_k(\zeta)$. We define the vector
$\boldx(\zeta,k,\boldd) = (a_1,\ldots,a_n)$ in $\mathbb{C}^n$, where $a_j = 1$ if $j$ is in $\mathbb{N}_k(\zeta)$ and $D_k(\varphi_j,\zeta) = \boldd$, while $a_j = 0$ otherwise. Let $\calM = \calM(\varphi_1,\ldots,\varphi_n)$ denote the linear span in $\mathbb{C}^n$ 
of all such vectors $\boldx(\zeta,k,\boldd)$. Clearly, condition (ii) in Theorem 5.13 is equivalent to
$$
(c_1,\ldots,c_n) \in \calM^\perp.
\leqno(\rm ii)'
$$
Put another way, if we define a linear transformation
$\mathbb{A}: \ \mathbb{C}^n \rarrow 
\calB(\calD_\beta)/\calK$ by
$$
\mathbb{A}(c_1,\ldots,c_n) = [c_1C_{\varphi_{1}}+ \cdots + c_n C_{\varphi_{n}}],
$$
where $[B]$ denotes the coset of the operator $B$, then the content of Theorem 5.13 is that ker $\mathbb{A} = \calM^\perp$. This yields the immediate:
}
\end{remark}

\begin{corollary}
Fix $\varphi_1,\ldots,\varphi_n$ in $\calS$. Then the vector subspace of $\calB(\calD_\beta)/\calK$ spanned by the cosets $[C_{\varphi_{1}}],\ldots,[C_{\varphi_{n}}]$
has the same dimension as the subspace $\calM(\varphi_1,\ldots,\varphi_n)$ of $\mathbb{C}^n$. In particular, these cosets are linearly independent if and only if $\calM(\varphi_1,\ldots,\varphi_n) = \mathbb{C}^n$.
\end{corollary}

For maps $\varphi$ in $\calS$ with order of contact uniformly two, the linear fractional self-maps of $D$ play a special role. Let us write $\calS(2)$ for the collection of those $\varphi$ in $\calS$ which have order of contact two at each point of $F(\varphi)$. Further, we denote by $\calL$ the collection of linear fractional self-maps $\varphi$ of $D$ which are not automorphisms but which have $\|\varphi\|_\infty =1$. We note that any linear fractional map $\psi$ is determined by its second-order data 
$D_2(\psi,z_0) = (\psi(z_0),\psi'(z_0),\psi''(z_0))$ at any point $z_0$ of analyticity. Suppose now that $\varphi$ is in $\calS(2)$ and $\zeta_0$ is a point in $F(\varphi)$. Let $\varphi_0$ be the unique linear fractional map with $D_2(\varphi_0,\zeta_0) = D_2(\varphi,\zeta_0)$. Since the curvature of the parametric curve $\varphi(e^{i\theta})$ is determined by second-order data, the circle $\varphi_0(\partial D)$ is necessarily the osculating circle of this curve at the point $\varphi(\zeta_0)$. Thus $\varphi_0$ maps $D$ to $D$ and lies in $\calL$. The following result was established by the second author on the weighted Bergman spaces \cite{Mo}; here we extend it to $H^2$.

\begin{corollary}
Let $\varphi$ be in $\calS(2)$ with $F(\varphi) = \{\zeta_1,\ldots,\zeta_r\}$. For $i = 1,\ldots,r$ let $\varphi_i$ be the unique linear fractional map with $D_2(\varphi_i,\zeta_i) = D_2(\varphi,\zeta_i)$. Then
$$
C_\varphi \equiv C_{\varphi_{1}} + \cdots + C_{\varphi_{r}} \ (\mbox{\rm mod } \calK),
$$
where the operators act on $H^2$.
\end{corollary}

{\em Proof.}
This is immediate from Theorem 5.13 applied to the linear combination
$$
C_\varphi - C_{\varphi_{1}}-C_{\varphi_{2}} - \cdots - C_{\varphi_{r}}.
$$
\hfill $\Box$ \bigskip

Let us write $\calV_2$ for the vector subspace of $\calB(\calD_\beta)/\calK$ spanned by the cosets $[C_\varphi]$ with $\varphi$ in $\calS(2)$. According to Corollary 5.16, $\{[C_\varphi]: \varphi \in \calL\}$ is a spanning subset of $\calV_2$. That it is linearly independent as well follows immediately from Theorem 5.13, which gives the following:

\begin{corollary}
$\{[C_\varphi]: \varphi \in \calL\}$ is a basis for $\calV_2$.
\end{corollary}

In Theorem 5.13, the proof that (ii) implies (i) uses the local compact difference Theorem 4.4 to locally ``match up'' certain pairs $C_{\varphi_{j}}$ and $C_{\varphi_{\ell}}$ by local data. As one changes the locality or local data, the pairs can change. However, the method suggests that with some minimality hypothesis, the coefficients in a linear equation
$$
c_1C_{\varphi_{1}}+  \cdots + c_n C_{\varphi_{n}} \equiv 0 \ (\mbox{\rm mod } \calK)
$$
should be integers up to a common scalar factor. The next result makes this precise.

\begin{corollary}
Suppose that $\varphi, \varphi_1, \ldots,\varphi_n$ lie in $\calS$ and assume that $[C_{\varphi_{1}}]$,
$\ldots$, $[C_{\varphi_{n}}]$ are linearly independent in $\calB(\calD_\beta)/\calK$. If
$$
C_\varphi \equiv c_1C_{\varphi_{1}}+ \cdots + c_n C_{\varphi_{n}} \ (\mbox{\rm mod } \calK),
$$
then $c_1,\ldots,c_n$ are rational numbers.
\end{corollary}

{\em Proof.}
Let us write $\psi_1 = \varphi$ and $\psi_k = \varphi_{k-1}$, $b_1 = -1$ and $b_k = c_{k-1}$ for $k = 2,3,\ldots,n+1$, so that the equation mod $\calK$ in the statement becomes
$$
b_1 C_{\psi_{1}} + b_2 C_{\psi_{2}} + \cdots + b_{n+1} C_{\psi_{n+1}} \equiv 0 \
(\mbox{\rm mod } \calK).
$$
Linear independence of $[C_{\psi_{2}}],\ldots,[C_{\psi_{n+1}}]$ implies that $b_1 = -1$ uniquely determines $b_2,\ldots,b_{n+1}$. Thus according to Remark 5.14, $\calM = \calM(\psi_1,\ldots,\psi_n)$ has codimension one in $\mathbb{C}^{n+1}$. From the vectors $\boldx(\zeta,k,\boldd)$ which span $\calM$, select a basis $\boldx_1,\ldots,\boldx_n$ for $\calM$; the coordinates of each $\boldx_i$ are zeros and ones. On applying the Gram-Schmidt process to $\{\boldx_1,\ldots,\boldx_n\}$, we obtain an orthonormal basis for $\calM$ of the form $\boldz_i/\|\boldz_i\|$, $i = 1,2,\ldots,n$, where each $\boldz_i$ has rational coefficients. Let $P$ denote the orthogonal projection of $\mathbb{C}^{n+1}$ onto $\calM$ and write $\bolde_1,\bolde_2,\ldots,\bolde_{n+1}$ for the standard basis vectors in $\mathbb{C}^{n+1}$. For some $j$, $(I-P)\bolde_j$ is non-zero, and clearly
$$
(I-P) \bolde_j = \bolde_j - \sum^n_{i=1} \frac{\langle \bolde_j,\boldz_i\rangle}{\|\boldz_i\|^2}
\ \boldz_i
$$
is a vector $(r_1,r_2,\ldots,r_{n+1})$ with obviously rational coefficients. Since $\calM^\perp$ is one-dimensional, $(b_1,\ldots,b_{n+1}) = c(r_1,\ldots,r_{n+1})$ for some complex number $c$. Since $cr_1 = b_1 = -1$, $c$ is rational so that $b_2,\ldots,b_{n+1}$ are rational as well.
\hfill $\Box$
\bigskip

Finally, we note that certain weighted composition operators enjoy a decomposition analogous to that in Corollary 5.16, with the sum replaced by a linear combination.

\begin{proposition}
Let $\varphi$ be in $\calS(2)$ with $F(\varphi) = \{\zeta_1,\ldots,\zeta_r\}$ and suppose that $\varphi_i$, $i = 1,\ldots,r$, are linear fractional maps related to $\varphi$ as in Corollary 5.16. If $w$ is a bounded measurable function on $\partial D$ which is continuous at each point of $F(\varphi)$, then
$$
M_w C_\varphi \equiv w(\zeta_1) C_{\varphi_{1}} + \cdots + w(\zeta_r)C_{\varphi_{r}} \ (\mbox{mod } \calK),
$$
where the operators are considered as mapping $H^2$ to $L^2$.
\end{proposition}

{\em Proof.} Choose pairwise disjoint arcs $I_1,\ldots,I_r$ with union $\partial\mathbb{D}$ and the interior of $I_j$ containing $\zeta_j$ for $j = 1,\ldots,r$. We write $\chi_{_{j}}$ for $\chi_{_{I_{j}}}$. We have
$$
M_w C_\varphi = \sum^r_{j=1} \left[M_{(w-w(\zeta_{j}))\chi_{_{j}}} C_\varphi +w(\zeta_j) M_{\chi_{_{j}}} C_\varphi\right].
$$
The first term in each summand on the right is compact by Theorem 3.1. Moreover the operators
$$
M_{\chi_{_{j}}}(C_\varphi-C_{\varphi_{j}}) \mbox{ and } M_{(1-\chi_{_{j}})} C_{\varphi_{j}}
$$
are compact, the first by Theorem 4.4 and Proposition 5.10(c) and the second by Theorem 3.1;  the proposition follows. \hfill $\Box$
\bigskip

Note that left-multiplying all operators in the proposition by the orthogonal projection $P$ of $L^2$ onto $H^2$ gives an alternate statement in which the multiplication operator $M_w$ is replaced by the corresponding Toeplitz operator $T_w = PM_w|H^2$, and all operators act from $H^2$ to $H^2$.

\section{Arc-Connectedness in $\calC(\calS_0)$}
\setcounter{equation}{0}

In \cite{SSu1}, Shapiro and Sundberg posed the following interesting question: What is the relationship between the conditions
\begin{itemize}
\item[(a)]
$C_\varphi-C_\psi$ is compact;
\item[(b)]
$C_\varphi$ and $C_\psi$ lie in the same component of the topological space of composition operators on $\calD_\beta$?
\end{itemize}
The second author and Toews \cite{MoT} proposed the scheme discussed in Section 2.6 and used it to give examples of composition operators satisfying (b) but not (a). Bourdon \cite{B} presented analogous examples for linear fractional maps, showing that for $\varphi$ and $\psi$ in the class  $\calL$ (see Section 5.5), (a) corresponds to equal second-order data at the (common) point of contact with $\partial D$ whereas equal first-order data characterizes (b). The second author [16] proved the corresponding result on the Bergman spaces for $\varphi$ and $\psi$ in the class $\calS_0 \cap \calS(2)$ where $\calS_0$ is as defined below. Here we extend this result to $H^2$ and to higher orders of contact. First we record a general sufficient condition for (b), proved in \cite{Mo} for the Bergman spaces.

\begin{proposition}
Let us write $\varphi_t = t\varphi + (1-t)\psi$, $0 \leq t \leq 1$, where
$\varphi$ and $\psi$ are analytic self-maps of $D$. If $\rho = \left|\frac{\varphi-\psi}{1-\overline{\varphi} \psi}\right|$ is essentially bounded away from one on $\partial D$, then there is a constant $B$ such that
$$
\|C_{\varphi_{s}}-C_{\varphi_{r}}\|_{\calB(H^2)} \leq 
B|s-r|, \quad 0 \leq r < s \leq 1.
$$
\end{proposition}

\noindent
The proposition follows immediately from Theorem 2.2 and Lemma 4.1 with $G$ taken to be the appropriate set of full measure.

We write $\calS_0$ for the collection of those $\varphi$ in $\calS$ such that $\|\chi_{_{\partial D\setminus U}}\varphi\|_\infty <1$
 for every open subset $U$ of $\partial D$ containing 
$F(\varphi)$. Such a $\varphi$ is allowed to make contact with $\partial D$ only at points of $F(\varphi)$. We also write $\calC(\calS_0)$ for the set of those $C_\varphi$ with $\varphi$ in $\calS_0$.

\begin{theorem}
Suppose that $\varphi$ and $\psi$ lie in $\calS_0$ with respective orders of contact $2m(\zeta)$ and $2n(\zeta)$ at each $\zeta$ in (respectively) $F(\varphi)$ and $F(\psi)$. Let $\varphi_t = t \varphi + (1-t)\psi$, $0 \leq t \leq 1$. Then if
$\beta \geq 1$, the following are equivalent:
\begin{itemize}
\item[(i)]
$C_\varphi$ and $C_\psi$ lie in the same component of the space $\calC(\calS_0)$ equipped with the norm topology of $\calB(\calD_\beta)$.
\item[(ii)]
There is a positive constant $B$ such that
$$
\|C_{\varphi_{s}}-C_{\varphi_{r}}\| \leq B|s-r|, \quad 0 \leq r < s \leq 1.
$$
\item[(iii)]
$F(\varphi) = F(\psi)$ (call this set $F$), and for each $\zeta$ in $F$, $2m(\zeta) = 2n(\zeta)$ and $D_{2m(\zeta)-1}(\varphi,\zeta) = D_{2m(\zeta)-1}(\psi,\zeta)$.
\end{itemize}
\end{theorem}

{\em Proof.} 
First consider the $H^2$ case. 
Suppose that (iii) holds and write $\rho = \left|\frac{\varphi-\psi}{1-\overline{\varphi}\psi}\right|$ as usual. If $\zeta$ is in $F$, it is clear from Proposition 5.10 that $\lim_{e^{i\theta}\rarrow \zeta} \rho(e^{i\theta}) <1$. Thus there is an open set $U$ in $\partial D$ containing $F$ such that $\|\chi_{_{U}} \rho \|_\infty <1$. Since $\varphi$ and $\psi$ lie in $\calS_0$, $\|\chi_{_{\partial D \setminus U}} \rho\|_\infty <1$ as well, and (ii) follows from Proposition 6.1.
For the Bergman space version of this implication, replace $U$ by an appropriate subset of $D$ and use the Bergman space analogue of Proposition 6.1 in \cite{Mo}, together with the following replacement for Proposition 5.10(b), a fact easily established by calculations in the proofs of Lemma 5.6 and Proposition 5.10: If $\varphi$ and $\psi$ have sufficient data and common order of contact $2m$ at $\zeta$ in $\partial D$, and if $D_{2m-1}(\varphi,\zeta) = D_{2m-1}(\psi,\zeta)$, then
$$
\limsup_{z \rarrow \zeta} \rho(z) < 1,
$$
the limit superior being taken unrestrictedly in $D$.
Clearly (ii) implies (i) on any $\calD_\beta$. 
Finally, we can use the topological argument of Theorem 2.4 in 
 \cite{M}, with Corollary 5.8 playing the role of Theorem 2.2 in \cite{M}, to show that (i) implies (iii). 
\hfill $\Box$ \bigskip

The last theorem and Corollary 5.12 show that the phenomenon (a) $\Rightarrow$ (b) (and not  conversely) persists more broadly. Note, however, that by passing from the class $\calS$ to $\calS_0$, we have eliminated maps $\varphi$ which touch the unit circle outside of $E(\varphi)$. Such points of contact are immaterial for condition (a); we do not know whether they matter for (b).

\end{document}